\documentclass[11pt, a4paper]{article}

\usepackage{array}
\usepackage{multirow}
\usepackage{arydshln}
\makeatletter
\newcommand{\thickhline}{%
    \noalign {\ifnum 0=`}\fi \hrule height 1pt
    \futurelet \reserved@a \@xhline
}
\newcolumntype{"}{@{\hskip\tabcolsep\vrule width 1pt\hskip\tabcolsep}}
\makeatother

\usepackage{graphics}
\usepackage{graphicx}
\usepackage{epsfig}
\usepackage{subfigure}
\usepackage{amssymb}
\usepackage{amsthm}
\usepackage{mathrsfs}
\usepackage{hhline}
\usepackage{longtable}
\usepackage{amsmath}
\usepackage{array}
\usepackage{float}
\usepackage{multirow}
\usepackage{cite}
\usepackage{enumerate}
\usepackage{xcolor}

\def\red{\color{red}}

\def\magenta{\color{magenta}}

\def\rsq{\hspace*{\fill}$\Box$\medskip}
\def\qed{\hspace*{\fill}$\blacksquare$\medskip}
\def\pf{\noindent{\bf Proof: }}
\def\Z{\mathbb Z}
\newtheoremstyle{de}
  {10pt}          
  {10pt}  
  {\rm}  
  {}
  {\bf}  
  {. }    
  { }    
  {}     
\theoremstyle{de}

\newtheorem{de}{Definition}[section]
\newtheorem{example}{Example}[section]
\newtheorem{rem}[de]{Remark}
\newtheorem{problem}{Problem}[section]
\newtheorem{question}{Question}[section]
\newtheoremstyle{theorem}
  {10pt}          
  {10pt}  
  {\it}  
  {}
  {\bf}  
  {. }    
  { }    
  {}     
\theoremstyle{theorem}
\newtheorem{theorem}{Theorem}[section]
\newtheorem{lemma}[theorem]{Lemma}
\newtheorem{proposition}[theorem]{Proposition}
\newtheorem{corollary}[theorem]{Corollary}

\numberwithin{equation}{section}
\def\Z{\mathbb{Z}}

\graphicspath{%
    {converted_graphics/}
    {/}
}

\begingroup\makeatletter\ifx\SetFigFont\undefined%
\gdef\SetFigFont#1#2#3#4#5{%
  \reset@font\fontsize{#1}{#2pt}%
  \fontfamily{#3}\fontseries{#4}\fontshape{#5}%
  \selectfont}%
\fi\endgroup%

\begin{document}
\baselineskip18truept
\normalsize
\begin{center}
{\mathversion{bold}\Large \bf Approaches Which Output Infinitely Many Graphs With Small Local Antimagic Chromatic Number}

\bigskip
{\large  Gee-Choon Lau$^{a,}$\footnote{Corresponding author.}, Jianxi Li$^b$, Ho-Kuen Ng$^c$, Wai-Chee Shiu{$^{d,e}$}}\\

\medskip

\emph{{$^a$}Faculty of Computer \& Mathematical Sciences,}\\
\emph{Universiti Teknologi MARA (Segamat Campus),}\\
\emph{85000, Johor, Malaysia.}\\
\emph{geeclau@yahoo.com}\\

\medskip
\emph{{$^b$}School of Mathematics and Statistics, Minnan Normal University,}\\
\emph{363000, Fujian, P.R. China}\\
\emph{ptjxli@hotmail.com}\\

\medskip

\emph{{$^c$}Department of Mathematics, San Jos\'{e} State University,}\\
\emph{San Jos\'e CA 95192 USA.}\\
\emph{ho-kuen.ng@sjsu.edu}\\

\medskip

\emph{{$^d$}Department of Mathematics, The Chinese University of Hong Kong,\\ Shatin, Hong Kong.}\\
\emph{{$^e$}College of Global Talents, Beijing Institute of Technology,\\ Zhuhai, China.}\\
\emph{wcshiu@associate.hkbu.edu.hk}\\

\end{center}
%

\begin{abstract}
An edge labeling of a connected graph $G = (V, E)$ is said to be local antimagic if it is a bijection $f:E \to\{1,\ldots ,|E|\}$ such that for any pair of adjacent vertices $x$ and $y$, $f^+(x)\not= f^+(y)$, where the induced vertex label $f^+(x)= \sum f(e)$, with $e$ ranging over all the edges incident to $x$.  The local antimagic chromatic number of $G$, denoted by $\chi_{la}(G)$, is the minimum number of distinct induced vertex labels over all local antimagic labelings of $G$. In this paper, we (i) give a sufficient condition for a graph with one pendant to have $\chi_{la}\ge 3$. A necessary and sufficient condition for a graph to have $\chi_{la}=2$ is then obtained; (ii) give a sufficient condition for every circulant graph of even order to have $\chi_{la} = 3$; (iii) construct infinitely many bipartite and tripartite graphs with $\chi_{la} = 3$ by transformation of cycles; (iv) apply transformation of cycles to obtain infinitely many one-point union of regular (possibly circulant) or bi-regular graphs with $\chi_{la} = 2,3$. The work of this paper suggests many open problems on the local antimagic chromatic number of bipartite and tripartite graphs. \\


\end{abstract}

\section{Introduction}\label{intro}
A connected graph $G = (V, E)$ is said to be local antimagic if it admits a local antimagic edge labeling, i.e., a bijection $f : E \to \{1,\ldots ,|E|\}$ such that the induced vertex labeling  $f^+ : V \to \mathbb Z$ given by $f^+(x) = \sum f(e)$ (with $e$ ranging over all the edges incident to $x$) has the property that any two adjacent vertices have distinct induced vertex labels.  The number of distinct induced vertex labels under $f$ is denoted by $c(f)$, and is called the {\it color number} of $f$. The {\it local antimagic chromatic number} of $G$, denoted by $\chi_{la}(G)$, is $\min\{c(f) : f\mbox{ is a local antimagic labeling of } G\}$~\cite{Arumugam}. In this paper, we (i) give a sufficient condition for a graph with one pendant to have $\chi_{la}\ge 3$. A necessary and sufficient condition for a graph to have $\chi_{la}=2$ is then obtained; (ii) give a sufficient condition for every circulant graph of even order to have $\chi_{la} = 3$; (iii) construct infinitely many bipartite and tripartite graphs with $\chi_{la} = 3$ by transformation of cycles; (iv) apply transformation of cycles to obtain infinitely many one-point union of regular (possibly circulant) or bi-regular graphs with $\chi_{la} = 2,3$.


\begin{lemma}\label{lem-2part}{\rm\cite{LauNgShiu-chila}} Let $G$ be a graph of size $q$. Suppose there is a local antimagic labeling of $G$ inducing a $2$-coloring of $G$ with colors $x$ and $y$, where $x<y$. Let $X$ and $Y$ be the numbers of vertices of colors $x$ and $y$, respectively. Then $G$ is a bipartite graph whose sizes of parts are $X$ and $Y$ with $X>Y$, and $xX=yY= \frac{q(q+1)}{2}$. \end{lemma}

\noindent Lemma~\ref{lem-2part} implies that
\begin{corollary}\label{cor-2part} Suppose $G$ is a bipartite graph of $q$ edges with bipartition $(V_1,V_2)$. If $\chi_{la}(G)=2$, then $|V_1|\ne |V_2|$ and ${q+1\choose2}$ is divisible by both $|V_1|$ and $|V_2|$.\end{corollary}

\noindent Note that the converse of this corollary does not hold. Consider the graph $G$ obtained from $P_7=u_1u_2\cdots u_7$ by adding edges $u_1u_4$ and $u_2u_5$. Thus, $G$ has bipartition
$V_1=\{u_1, u_3, u_5, u_7\}$ and $V_2=\{u_2, u_4, u_6\}$ and $q=8$ satisfying the conclusion of Corollary~\ref{cor-2part}. If $\chi_{la}(G)=2$, then vertex $u_7$ must have label 9, which is impossible.
\begin{theorem}\label{thm-2parteven} Suppose $G$ is a bipartite graph of even size $q$ and contains one pendant, then $\chi_{la}(G)\ge 3$.
\end{theorem}
\pf Suppose $\chi_{la}(G)= 2$. Let $(V_1, V_2)$ be a bipartition of $G$. Then the pendant edge must be labeled by $q$. By Lemma~\ref{lem-2part}, $q|V_1|=\frac{q(q+1)}{2}$, where $V_1$ contains the pendant  vertex. It is impossible.
\rsq

\noindent The following theorem gives a necessary and sufficient condition for a graph $G$ to have $\chi_{la}(G)=2$.

\begin{theorem}\label{thm-2part} A graph $G$ of size $q$ has $\chi_{la}(G)=2$ if and only if $G$ is bipartite with bipartition $(V_1, V_2)$ such that (i) $G$ has at most one pendant, and (ii) $G$ admits a local antimagic labeling with every vertex in $V_1$ has label ${m+1 \choose 2}/|V_1|$ and every vertex in $V_2$ has label ${m+1 \choose 2}/|V_2|$ that are distinct integers.   \end{theorem}

\pf If $\chi_{la}(G)=2$, Lemma~\ref{lem-2part} implies that $G$ is bipartite with every vertex in $V_1$ has label ${m+1 \choose 2}/|V_1|$ and every vertex in $V_2$ has label ${m+1 \choose 2}/|V_2|$ that are distinct integers. Clearly, $G$ cannot have at least 2 pendants, otherwise, $\chi_{la}(G)\ge 3$. Conversely, if $G$ satisfies the given conditions, then $\chi_{la}(G)\le 2$. Since $\chi_{la}(G)\ge \chi(G)=2$, the equality holds.
\rsq

\begin{lemma}\label{lem-reg}{\rm\cite{LauNgShiu-chila}} Suppose $G$ is a $d$-regular graph of size $q$. If $f$ is a local antimagic labeling of $G$, then $g = q + 1 - f$ is also a local antimagic labeling of $G$ with $c(f)= c(g)$.  Moreover, suppose $c(f)=\chi_{la}(G)$ and if $f(uv)=1$ or $f(uv)=q$, then $\chi_{la}(G-uv)\le \chi_{la}(G)$. \end{lemma}

\begin{lemma}\label{lem-nonreg}{\rm\cite{LauNgShiu-chila}} Suppose $G$ is a graph of size $q$ and $f$ is a local antimagic labeling of $G$. For any $x,y\in V(G)$, if\\
(i) $f^+(x) = f^+(y)$ implies that $\deg(x)=\deg(y)$, and\\
(ii) $f^+(x) \ne f^+(y)$ implies that $(q+1)(\deg(x)-\deg(y)) \ne f^+(x) - f^+(y)$,\\ then $g = q + 1 - f$ is also a local antimagic labeling of $G$ with $c(f)= c(g)$. \end{lemma}

\noindent  For $t\ge 2$, consider the following conditions for a graph $G$:

\begin{enumerate}[(i)]
  \item $\chi_{la}(G)=t$ and $f$ is a local antimagic labeling of $G$ that induces a $t$-independent partition $\cup^t_{i=1} V_i$ of $V(G)$.
  \item For each $x\in V_k$, $1\le k\le t$, $deg(x)=d_k$ satisfying $f^+(x) - d_a \ne f^+(y) - d_b$, where $x\in V_a$ and $y\in V_b$ for $1\le a\ne b\le t$.
\end{enumerate}

\begin{lemma}\label{lem-G-e}{\rm\cite{LauNgShiu-chila}} Let $H$ be obtained from $G$ with an edge $e$ deleted. If $G$ satisfies Conditions (i) and (ii) and $f(e)=1$, then $\chi(H)\le \chi_{la}(H)\le t$. \end{lemma}

\section{Circulant Graphs}\label{sec-cir}

Suppose $A$ is an additive (Abelian) group and $S$ is a set of generators of $A$. The Cayley (simple undirected) graph $G=\Gamma(A, S)$ associated with $(A,S)$ is defined as follows:
\begin{enumerate}[1.]
\item The vertex set of $G$ is $A$.
\item The edge set of $G$ is $\{uv \; |\; u-v\in S \mbox{ or }v-u\in S\}$.
\end{enumerate}

\noindent For each $a\in S$, we denote $\Gamma_a=\Gamma(A, \{a\})$. Then $\Gamma(A, S)=\bigcup\limits_{a\in S} \Gamma_a$. Note that, $\Gamma_a=\Gamma_{-a}$.

\noindent In the section, we focus on some special circulant graphs which are Cayley graphs.
So, we fix $A=\Z_m=\{0,1,\dots, m-1\}$, where $m\ge 3$. Let $D=\{d\;|\;1\le d< \lceil m/2\rceil, (d,m)=1\}$, where $(a,b)$ denotes the g.c.d. of $a$ and $b$. Suppose $S=\{a_0, a_1,\dots, a_t\}\subseteq D$, where $t\ge 0$. If we rename the vertex $v_i$ by $v_{a_0^{-1}i}$, then $\Gamma_{a_i}\cong \Gamma_{a_0^{-1}a_i}$ for all $i\in \Z_m$. So, without loss of generality, we may assume that $a_0=1$.
The Cayley graph $\Gamma(\Z_m, S)$ is a circulant graph and is denoted by
$C_m(a_0,a_1, a_2, \ldots, a_t)$. Thus $\Gamma(\Z_m, \{1\})=C_m$, the $m$-cycle and $C_m(a_0,a_1, a_2, \ldots, a_t)=\bigcup\limits_{i=0}^t \Gamma_{a_i}$. For each $a_i$, $\Gamma_{a_i}$ is the $m$-cycle $(0, a_i, 2a_i, \dots, (m-1)a_i)$, $0\le i\le t$. Without loss of generality, we may assume that $1<a_1<\cdots <a_t$.

\begin{lemma}\label{lem-iso}
Suppose $(a,n)=1$ and $ab\equiv 1\pmod{n}$, then $C_n(1,a) \cong C_n(1,b)$.
\end{lemma}
\pf Let $g: \Z_n \to \Z_n$ be defined by $g(i)=ib$. Clearly $g$ is an automorphism, since $(b,n)=1$. For any vertex $j$, the neighborhood of $j$ is $N(j)=\{j-1, j+1, j-a, j+a\}$. Then $g$ maps $N(j)$ to $\{(j-1)b, (j+1)b, (j-a)b, (j+a)b\}=\{jb-b, jb+b, jb-ab, jb+ab\}=\{jb-b, jb+b, jb-1, jb+1\}$ which is the neighborhood of $jb$ in $C_n(b)$. Thus $g$ induces an isomorphism from $C_n(1,a)$ onto $C_n(1,b)$.
\rsq

\noindent Since $\Gamma_b=\Gamma_{-b}$, we have
\begin{corollary} Suppose $(a,n)=1$ and $ab\equiv -1\pmod{n}$, then $C_n(1,a) \cong C_n(1,b)$.\end{corollary}

\noindent For convenience, we also use $v_0, v_1, \dots, v_{m-1}$ instead of the $C_m$'s vertices $0, 1, \dots, m-1$, respectively.

\begin{example}\label{ex-C16} Take $n=16$, then we have $D=\{1,3,5,7\}$. All $C_{16}(1,a_1,\dots, a_t)$ are $C_{16}$, $C_{16}(1,3)$, $C_{16}(1,5)$, $C_{16}(1,7)$, $C_{16}(1,3,5)$, $ C_{16}(1,3,7)$, $C_{16}(1,5,7)$, $C_{16}(1,3,5,7)$.

\noindent By Lemma~\ref{lem-iso}, $C_{16}(1,3)$ is isomorphic to $C_{16}(1,11)=C_{16}(1,5)$. Define $\phi:\Z_{16}\to \Z_{16}$ by $\phi(i)=5i$ and $\psi:\Z_{16}\to \Z_{16}$ by $\psi(i)=3i$. It is easy to check that $\psi$ and $\phi$ induce isomorphisms from $C_{16}(1,3,5)$ to $C_{16}(1,3,7)$ and $C_{16}(1,5,7)$, respectively.

\noindent However, $C_{16}(1,3) \not\cong C_{16}(1,7)$. We may consider the spectra of $C_{16}(1,3)$ and $C_{16}(1,7)$. One may find the formula of the spectrum of a circulant matrix from \cite{ShiuMaFang}. Let $\rho_1(x)=x+x^{15}$, $\rho_3(x)=x^3+x^{13}$ and $\rho_7(x)=x^7+x^9$. Let $\zeta=e^{2\pi i/16}$, then the spectra of $C_{16}(3)$ and $C_{16}(7)$ are
\[\{\rho_1(\zeta^j)+\rho_3(\zeta^j)\;|\; 0\le j\le 15\}\mbox{ and } \{\rho_1(\zeta^j)+\rho_7(\zeta^j)\;|\; 0\le j\le 15\},\]
respectively. Since $\rho_k(\zeta^r)=\rho_k(\zeta^{-r})$ for $1\le r\le 15$, we only list $\rho_1(\zeta^j)+\rho_3(\zeta^j)$ and $\rho_1(\zeta^j)+\rho_7(\zeta^j)$ for $0\le j\le 8$. By direct computation we have
\[\begin{array}{r|l}
j& \rho_1(\zeta^j)+\rho_3(\zeta^j)\\\hline
0 & 4\\
1 & \zeta+\zeta^{3}+\zeta^{-3}+\zeta^{-1}\\
2 & \zeta^2+\zeta^{6}+\zeta^{-6}+\zeta^{-2}\\
3 & \zeta^3+\zeta^{7}+\zeta^{-7}+\zeta^{-3}\\
4 & 0\\
5 & \zeta+\zeta^5+\zeta^{-5}+\zeta^{-1}\\
6 & \zeta^2+\zeta^{6}+\zeta^{-6}+\zeta^{-2}\\
7 & \zeta^5+\zeta^{7}+\zeta^{-7}+\zeta^{-5}\\
8 & -4
\end{array}\qquad
\begin{array}{r|l}
j& \rho_1(\zeta^j)+\rho_7(\zeta^j)\\\hline
0 & 4\\
1 & \zeta+\zeta^{7}+\zeta^{-7}+\zeta^{-1}\\
2 & 2(\zeta^2+\zeta^{-2})\\
3 & \zeta^3+\zeta^{5}+\zeta^{-5}+\zeta^{-3}\\
4 & 0\\
5 & \zeta^3+\zeta^{5}+\zeta^{-5}+\zeta^{-3}\\
6 & 2(\zeta^6+\zeta^{-6})\\
7 & \zeta+\zeta^{7}+\zeta^{-7}+\zeta^{-1}\\
8 & -4
\end{array}
\]
\noindent Clearly, they are not the same. So $C_{16}(1,3)$ and $C_{16}(1,7)$ are not isomorphic.
\qed
\end{example}

\noindent In~\cite{Arumugam}, it was shown that $\chi_{la}(C_m) = 3$ for all $m\ge 3$. Throughout this paper, we shall refer to the following local antimagic labeling  $f$ for a cycle $C_m$, denoted {\it $C$-labeling}, whenever necessary:

$$f(v_jv_{j+1})=\begin{cases} (j+2)/2 & \mbox{ for even } j, \\ m - (j-1)/2 & \mbox{ for odd } j\end{cases}$$ so that $f^+(v_0) = \lfloor \frac{m}{2}\rfloor+2$, $f^+(v_{j}) = m+1$ for odd $j$, and $f^+(v_j) = m+2$ for even $j > 0$.

\medskip

\noindent Consider an $m$-cycle $(0, a, 2a, \dots, (m-1)a) =\Gamma_{a}$. Now for each $i\ge 0$, we label the $m$-cycle $\Gamma_{a}$ according to the $C$-labeling translating by $im$, that is, $v_{ja}v_{(j+1)a}$ is labeled by $f(v_jv_{j+1})+im$. Now $\Gamma_{a}$ is labeled by integers in $[im+1, (i+1)m]$. We denote this labeling by $f_i$. Note that $f_0=f$. Thus, $f_i^+(v_0) = \lfloor \frac{m}{2}\rfloor+2+2im$, $f_i^+(v_{ja}) = m+1+2im$ for odd $j$, and $f_i^+(v_{ja}) = m+2+2im$ for even $j > 0$.

\noindent Now we consider the graph $C_m(a_0, a_1, a_2, \ldots, a_t)$. We label each $\Gamma_{a_i}$ by $f_i$.
Combining the labelings $f_0, f_1,\dots f_t$ we have a labeling $g$ for the whole graph $C_m(a_0, a_1, a_2, \ldots, a_t)$ given by
\[ g(v_{ja_i}v_{(j+1)a_i})=f(v_jv_{j+1})+im,\quad 0\le i\le t,\ 0\le j\le m-1.\]

\begin{example} Let $m=9$ so that $D=\{1,2,4\}$. Now, $\Gamma_1=(0,1,2,3,4,5,6,7,8)$, the $9$-cycle in natural order, $\Gamma_2=(0,2,4,6,8,1,3,5,7)$ and $\Gamma_4=(0,4,8,3,7,2,6,1,5)$.

\noindent Suppose $S=\{1,2\}$. According to the labelings defined above, we have the following vertex labelings:
\[\begin{array}{c||*{9}{c|}}
& 0 & 1 & 2 & 3 & 4 & 5 & 6 & 7 & 8\\\hline
f^+_0 & 6 & 10 & 11 & 10 & 11 & 10 & 11 & 10 & 11\\
f^+_1 & 24 & 28 & 28 & 29 & 29 & 28 & 28 & 29 & 29\\\hline
g & 30 & 38 & 39 & 39 & 40 & 38 & 39 & 39 & 40
\end{array}\]

\noindent Suppose $S=D$.  According to the labelings defined above, we have the following vertex labelings:
\[\begin{array}{c||*{9}{c|}}
& 0 & 1 & 2 & 3 & 4 & 5 & 6 & 7 & 8\\\hline
f^+_0 & 6 & 10 & 11 & 10 & 11 & 10 & 11 & 10 & 11\\
f^+_1 & 24 & 28 & 28 & 29 & 29 & 28 & 28 & 29 & 29\\
f^+_2 & 42 & 46 & 46 & 46 & 46 & 47 & 47 & 47 & 47\\\hline
g & 72 & 84 & 85 & 85 & 86 & 85 & 86 & 86 & 87
\end{array}\]
\noindent Here both $g$'s are not local antimagic labelings.
\qed
\end{example}

\noindent Now let us consider even $m$, i.e., $m=2n$. In this case, all $a_i$ are odd. So $f_{i}^+(v_0) = n+2+4in$, $f_{i}^+(v_{j}) = 2n+1+4in$ for odd $j$, and $f_i^+(v_j) = 2n+2+4in$ for even $j > 0$. Hence $g^+(v_0) = (t+1)(2nt+n+2)$, $g^+(v_{j}) = (t+1)(2nt+2n+1)$ for odd $j$, and $g^+(v_j) = (t+1)(2nt+2n+2)$ for even $j > 0$. Thus, together with the contrapositive of Lemma~\ref{lem-2part}, we have

\begin{theorem}\label{thm-circulant} For $1<a_1<\cdots < a_t<n$ and $(a_j, 2n)=1$, $1\le j\le t$, $\chi_{la}(C_{2n}(1,a_1, \ldots, a_t))=3$.
\end{theorem}

\begin{corollary} For each even $r\ge2$, there are infinitely many $r$-regular graphs of even order with local antimagic chromatic number equal $3$ but chromatic number equal $2$. \end{corollary}

\noindent Note that each $C_{2n}(1,a_1, a_2, \ldots, a_t)$ is edge transitive. Thus, all the graphs obtained by deleting one edge are isomorphic. By Lemmas~\ref{lem-2part} and~\ref{lem-reg}, the following is obvious.

\begin{corollary} The graph $C_{2n}(1,a_1, a_2, \ldots, a_t)$ with an edge deleted has local antimagic chromatic number equal $3$, where $1<a_1<\cdots < a_t<n$ and $(a_j, 2n)=1$ for $1\le j\le t$. \end{corollary}

\begin{figure}[H]
\begin{center}
\scalebox{0.4}{
\begin{picture}(0,0)%
\includegraphics{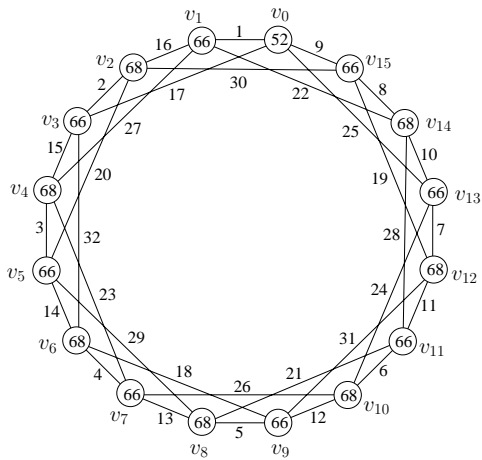}%
\end{picture}%

\begin{picture}(6989,7174)(569,-6055)
\put(6159,143){\makebox(0,0)[lb]{\smash{{\SetFigFont{20}{24.0}{\familydefault}{\mddefault}{\updefault}{\color[rgb]{0,0,0}$v_{15}$}%
}}}}
\put(7092,-843){\makebox(0,0)[lb]{\smash{{\SetFigFont{20}{24.0}{\familydefault}{\mddefault}{\updefault}{\color[rgb]{0,0,0}$v_{14}$}%
}}}}
\put(7491,-3224){\makebox(0,0)[lb]{\smash{{\SetFigFont{20}{24.0}{\familydefault}{\mddefault}{\updefault}{\color[rgb]{0,0,0}$v_{12}$}%
}}}}
\put(6127,-5280){\makebox(0,0)[lb]{\smash{{\SetFigFont{20}{24.0}{\familydefault}{\mddefault}{\updefault}{\color[rgb]{0,0,0}$v_{10}$}%
}}}}
\put(7543,-1976){\makebox(0,0)[lb]{\smash{{\SetFigFont{20}{24.0}{\familydefault}{\mddefault}{\updefault}{\color[rgb]{0,0,0}$v_{13}$}%
}}}}
\put(3344,845){\makebox(0,0)[lb]{\smash{{\SetFigFont{20}{24.0}{\familydefault}{\mddefault}{\updefault}{\color[rgb]{0,0,0}$v_1$}%
}}}}
\put(4689,876){\makebox(0,0)[lb]{\smash{{\SetFigFont{20}{24.0}{\familydefault}{\mddefault}{\updefault}{\color[rgb]{0,0,0}$v_0$}%
}}}}
\put(3428,-5943){\makebox(0,0)[lb]{\smash{{\SetFigFont{20}{24.0}{\familydefault}{\mddefault}{\updefault}{\color[rgb]{0,0,0}$v_8$}%
}}}}
\put(4693,-5955){\makebox(0,0)[lb]{\smash{{\SetFigFont{20}{24.0}{\familydefault}{\mddefault}{\updefault}{\color[rgb]{0,0,0}$v_9$}%
}}}}
\put(2186,-5429){\makebox(0,0)[lb]{\smash{{\SetFigFont{20}{24.0}{\familydefault}{\mddefault}{\updefault}{\color[rgb]{0,0,0}$v_7$}%
}}}}
\put(639,-1880){\makebox(0,0)[lb]{\smash{{\SetFigFont{20}{24.0}{\familydefault}{\mddefault}{\updefault}{\color[rgb]{0,0,0}$v_4$}%
}}}}
\put(584,-3190){\makebox(0,0)[lb]{\smash{{\SetFigFont{20}{24.0}{\familydefault}{\mddefault}{\updefault}{\color[rgb]{0,0,0}$v_5$}%
}}}}
\put(1078,-4310){\makebox(0,0)[lb]{\smash{{\SetFigFont{20}{24.0}{\familydefault}{\mddefault}{\updefault}{\color[rgb]{0,0,0}$v_6$}%
}}}}
\put(1122,-794){\makebox(0,0)[lb]{\smash{{\SetFigFont{20}{24.0}{\familydefault}{\mddefault}{\updefault}{\color[rgb]{0,0,0}$v_3$}%
}}}}
\put(1968,147){\makebox(0,0)[lb]{\smash{{\SetFigFont{20}{24.0}{\familydefault}{\mddefault}{\updefault}{\color[rgb]{0,0,0}$v_2$}%
}}}}
\put(6988,-4367){\makebox(0,0)[lb]{\smash{{\SetFigFont{20}{24.0}{\familydefault}{\mddefault}{\updefault}{\color[rgb]{0,0,0}$v_{11}$}%
}}}}
\end{picture}
}
\caption{A local antimagic labeling for $C_{16}(1,3)$ with color number $3$.}\label{fig-C16-3}
\end{center}
\end{figure}

\begin{example} By Theorem~\ref{thm-circulant} the local antimagic chromatic number of each graph listed in Example~\ref{ex-C16} is $3$.  As an example to illustrate the proof of Theorem~\ref{thm-circulant}, let us consider $C_{16}(1,3)$. Now $\Gamma_3=(0,3,6,9,12,15,2,5,8,11,14,1,4,7,10,13)$. The corresponding local antimagic labeling of $C_{16}(1,3)$ is given in Figure~\ref{fig-C16-3}. \qed
\end{example}
%

\noindent Let us use labeling matrices $M_f(\Gamma_1), M_{f_1}(\Gamma_3), M_{f_1}(\Gamma_7)$ to represent the labeling for $\Gamma_1, \Gamma_3, \Gamma_7$ according to the labeling $f$ and $f_{i}$ defined above, respectively.
\[M_f(\Gamma_1)= \fontsize{9}{11}\selectfont
\begin{array}{r|*{16}{p{2mm}}|r}
& 0 & 1 & 2 & 3 & 4 & 5 & 6 & 7 & 8 & 9 & 10 & 11 & 12 & 13 & 14 & 15 & \mbox{Sum}\\\hline
0 & * & 1 & * & * & * & * & * & * & * & * & * & * & * & * & * & 9 & 10\\
1 & 1 & * & 16 & * & * & * & * & * & * & * & * & * & * & * & * & * & 17\\
2 & * & 16 & * & 2 & * & * & * & * & * & * & * & * & * & * & * & * & 18\\
3 & * & * & 2 & * & 15 & * & * & * & * & * & * & * & * & * & * & * & 17\\
4 & * & * & * & 15 & * & 3 & * & * & * & * & * & * & * & * & * & * & 18\\
5 & * & * & * & * & 3 & * & 14 & * & * & * & * & * & * & * & * & * & 17\\
6 & * & * & * & * & * & 14 & * & 4 & * & * & * & * & * & * & * & * & 18\\
7 & * & * & * & * & * & * & 4 & * & 13 & * & * & * & * & * & * & * & 17\\
8 & * & * & * & * & * & * & * & 13 & * & 5 & * & * & * & * & * & * & 18\\
9 & * & * & * & * & * & * & * & * & 5 & * & 12 & * & * & * & * & * & 17\\
10 & * & * & * & * & * & * & * & * & * & 12 & * & 6 & * & * & * & * & 18\\
11 & * & * & * & * & * & * & * & * & * & * & 6 & * & 11 & * & * & * & 17\\
12 & * & * & * & * & * & * & * & * & * & * & * & 11 & * & 7 & * & * & 18\\
13 & * & * & * & * & * & * & * & * & * & * & * & * & 7 & * & 10 & * & 17\\
14 & * & * & * & * & * & * & * & * & * & * & * & * & * & 10 & * & 8 & 18\\
15 & 9 & * & * & * & * & * & * & * & * & * & * & * & * & * & 8 & * & 17\\
\end{array}
\]
\[M_{f_1}(\Gamma_3)=\fontsize{9}{11}\selectfont
\begin{array}{r|*{16}{p{2mm}}|r}
& 0 & 1 & 2 & 3 & 4 & 5 & 6 & 7 & 8 & 9 & 10 & 11 & 12 & 13 & 14 & 15 & \mbox{Sum}\\\hline
0 & * & * & * & 17 & * & * & * & * & * & * & * & * & * & 25 & * & * & 42\\
1 & * & * & * & * & 27 & * & * & * & * & * & * & * & * & * & 22 & * & 49\\
2 & * & * & * & * & * & 20 & * & * & * & * & * & * & * & * & * & 30 & 50\\
3 & 17 & * & * & * & * & * & 32 & * & * & * & * & * & * & * & * & * & 49\\
4 & * & 27 & * & * & * & * & * & 23 & * & * & * & * & * & * & * & * & 50\\
5 & * & * & 20 & * & * & * & * & * & 29 & * & * & * & * & * & * & * & 49\\
6 & * & * & * & 32 & * & * & * & * & * & 18 & * & * & * & * & * & * & 50\\
7 & * & * & * & * & 23 & * & * & * & * & * & 26 & * & * & * & * & * & 49\\
8 & * & * & * & * & * & 29 & * & * & * & * & * & 21 & * & * & * & * & 50\\
9 & * & * & * & * & * & * & 18 & * & * & * & * & * & 31 & * & * & * & 49\\
10 & * & * & * & * & * & * & * & 26 & * & * & * & * & * & 24 & * & * & 50\\
11 & * & * & * & * & * & * & * & * & 21 & * & * & * & * & * & 28 & * & 49\\
12 & * & * & * & * & * & * & * & * & * & 31 & * & * & * & * & * & 19 & 50\\
13 & 25 & * & * & * & * & * & * & * & * & * & 24 & * & * & * & * & * & 49\\
14 & * & 22 & * & * & * & * & * & * & * & * & * & 28 & * & * & * & * & 50\\
15 & * & * & 30 & * & * & * & * & * & * & * & * & * & 19 & * & * & * & 49\\
\end{array}
\]
\[M_{f_1}(\Gamma_7)=\fontsize{9}{11}\selectfont
\begin{array}{r|*{16}{p{2mm}}|r}
& 0 & 1 & 2 & 3 & 4 & 5 & 6 & 7 & 8 & 9 & 10 & 11 & 12 & 13 & 14 & 15 & \mbox{Sum}\\\hline
0 & * & * & * & * & * & * & * & 17 & * &  25 & * & * & * & * & * & * & 42\\
1 & * & * & * & * & * & * & * & * & 29 & * & 20 & * & * & * & * & * & 49\\
2 & * & * & * & * & * & * & * & * & * & 24 & * & 26 & * & * & * & * & 50\\
3 & * & * & * & * & * & * & * & * & * & * & 30 & * & 19 & * & * & * & 49\\
4 & * & * & * & * & * & * & * & * & * & * & * & 23 & * & 27 & * & * & 50\\
5 & * & * & * & * & * & * & * & * & * & * & * & * & 31 & * & 18 & * & 49\\
6 & * & * & * & * & * & * & * & * & * & * & * & * & * & 22 & * & 28 & 50\\
7 & 17 & * & * & * & * & * & * & * & * & * & * & * & * & * & 32 &  * & 49\\
8 & * & 29 & * & * & * & * & * & * & * & * & * & * & * & * & * & 21 & 50\\
9 & 25 & * & 24 & * & * & * & * & * & * & * & * & * & * & * & * & * & 49\\
10& * & 20 & * & 30 & * & * & * & * & * & * & * & * & * & * & * & * & 50\\
11 & * & * & 26 & * & 23 & * &  * & * & * & * & * & * & * & * & * & * & 49\\
12 & * & * & * & 19 & * & 31 & * & * & * & * & * & * & * & * & * & * & 50\\
13 & * & * & * & * & 27 & * & 22 & * & * & * & * & * & * & * & * & * & 49\\
14 & * & * & * & * & * & 18 & * & 32 & * & * & * & * & * & * & * & * & 50\\
15 & * & * & * & * & * & * & 28 & * & 21 & * & * & * & * & * & * & * & 49\\
\end{array}
\]
\noindent Note that $M_{f_{i}}(\Gamma_j)$ is the matrix obtained from $M_{f}(\Gamma_j)$ by adding all numerical entries by $im$. Here, we do not show the labeling matrices for the labelings $f_2$ and $f_3$.

\noindent According to the labeling $g$ for $C_{16}(1,3)$ defined above, the labeling matrix is $M_g( C_{16}(1,3))=M_f(\Gamma_1)+M_{f_1}(\Gamma_3)$, where $*$ is treated as 0:
\[M_g( C_{16}(1,3))=\fontsize{9}{11}\selectfont
\begin{array}{r|*{16}{p{2mm}}|r}
& 0 & 1 & 2 & 3 & 4 & 5 & 6 & 7 & 8 & 9 & 10 & 11 & 12 & 13 & 14 & 15 & \mbox{Sum}\\\hline
0 & * & 1 & * & 17 & * & * & * & * & * & * & * & * & * & 25 & * & 9 & 52\\
1 & 1 & * & 16 & * & 27 & * & * & * & * & * & * & * & * & * & 22 & * & 66\\
2 & * & 16 & * & 2 & * & 20 & * & * & * & * & * & * & * & * & * & 30 & 68\\
3 & 17 & * & 2 & * & 15 & * & 32 & * & * & * & * & * & * & * & * & * & 66\\
4 & * & 27 & * & 15 & * & 3 & * & 23 & * & * & * & * & * & * & * & * & 68\\
5 & * & * & 20 & * & 3 & * & 14 & * & 29 & * & * & * & * & * & * & * & 66\\
6 & * & * & * & 32 & * & 14 & * & 4 & * & 18 & * & * & * & * & * & * & 68\\
7 & * & * & * & * & 23 & * & 4 & * & 13 & * & 26 & * & * & * & * & * & 66\\
8 & * & * & * & * & * & 29 & * & 13 & * & 5 & * & 21 & * & * & * & * & 68\\
9 & * & * & * & * & * & * & 18 & * & 5 & * & 12 & * & 31 & * & * & * & 66\\
10 & * & * & * & * & * & * & * & 26 & * & 12 & * & 6 & * & 24 & * & * & 68\\
11 & * & * & * & * & * & * & * & * & 21 & * & 6 & * & 11 & * & 28 & * & 66\\
12 & * & * & * & * & * & * & * & * & * & 31 & * & 11 & * & 7 & * & 19 & 68\\
13 & 25 & * & * & * & * & * & * & * & * & * & 24 & * & 7 & * & 10 & * & 66\\
14 & * & 22 & * & * & * & * & * & * & * & * & * & 28 & * & 10 & * & 8 & 68\\
15 & 9 & * & 30 & * & * & * & * & * & * & * & * & * & 19 & * & 8 & * & 66\\
\end{array}
\]
Also $M_g( C_{16}(1,7))=M_f(\Gamma_1)+M_{f_1}(\Gamma_7)$ which is
\[\fontsize{9}{11}\selectfont
\begin{array}{r|*{16}{p{2mm}}|r}
& 0 & 1 & 2 & 3 & 4 & 5 & 6 & 7 & 8 & 9 & 10 & 11 & 12 & 13 & 14 & 15 & \mbox{Sum}\\\hline
0 & * & 1 & * & * & * & * & * & 17 & * &  25 & * & * & * & * & * & 9 & 52\\
1 & 1 & * & 16 & * & * & * & * & * & 29 & * & 20 & * & * & * & * & * & 66\\
2 & * & 16 & * & 2 & * & * & * & * & * & 24 & * & 26 & * & * & * & * & 68\\
3 & * & * & 2 & * & 15 & * & * & * & * & * & 30 & * & 19 & * & * & * & 66\\
4 & * & * & * & 15 & * & 3 & * & * & * & * & * & 23 & * & 27 & * & * & 68\\
5 & * & * & * & * & 3 & * & 14 & * & * & * & * & * & 31 & * & 18 & * & 66\\
6 & * & * & * & * & * & 14 & * & 4 & * & * & * & * & * & 22 & * & 28 & 68\\
7 & 17 & * & * & * & * & * & 4 & * & 13 & * & * & * & * & * & 32 &  * & 66\\
8 & * & 29 & * & * & * & * & * & 13 & * & 5 & * & * & * & * & * & 21 & 68\\
9 & 25 & * & 24 & * & * & * & * & * & 5 & * & 12 & * & * & * & * & * & 66\\
10& * & 20 & * & 30 & * & * & * & * & * & 12 & * & 6 & * & * & * & * & 68\\
11 & * & * & 26 & * & 23 & * &  * & * & * & * & 6 & * & 11 & * & * & * & 66\\
12 & * & * & * & 19 & * & 31 & * & * & * & * & * & 11 & * & 7 & * & * & 68\\
13 & * & * & * & * & 27 & * & 22 & * & * & * & * & * & 7 & * & 10 & * & 66\\
14 & * & * & * & * & * & 18 & * & 32 & * & * & * & * & * & 10 & * & 8 & 68\\
15 & 9 & * & * & * & * & * & 28 & * & 21 & * & * & * & * & * & 8 & * & 66\\
\end{array}
\]

\noindent Similarly, we have\\ $M_g(C_{16}(1,3,7))=M_f(\Gamma_1)+M_{f_1}(\Gamma_3)+M_{f_2}(\Gamma_7)$. One may choose $M_f(\Gamma_1)+M_{f_2}(\Gamma_3)+M_{f_1}(\Gamma_7)$ or other combination.\\ $M_g(C_{16}(1,3,5,7))=M_f(\Gamma_1)+M_{f_1}(\Gamma_3)+M_{f_2}(\Gamma_5)+M_{f_3}(\Gamma_7)$.  We do not show the matrices here. \rsq

\begin{rem} Suppose $a$ is not a generator of a finite Abelian group $A$. We can still define a graph $\Gamma_a$ with vertex set $A$ and the edge set $\{uv\;|\; u-v=a\mbox{ or }v-u=a\}$. The graph $\Gamma_a$ is a disjoint union of cycles. Precisely, let $H$ be the cyclic group generated by $a$ of order $n$. Then $A$ is the union of its cosets, namely $\bigcup\limits_{i=1}^s (a_i+H)$ for some $a_i\in A$, where $s$ is the index of the subgroup $H$ of $A$. Each coset $a_i+H=\{a_i, a_i+a, \dots, a_i+(n-1)a\}$ corresponds an $n$-cycle $(a_i, a_i+a, \dots, a_i+(n-1)a)=\Lambda_{a_i}$. Then $\Gamma_a=\sum\limits_{i=1}^s \Lambda_{a_i}$.
\end{rem}

\section{Transformation of Cycles}\label{cycles}

 We now give ways to transform an $n$-cycle $C_n=v_0v_1v_2\cdots v_{n-1}v_0$ that has a $C$-labeling $f$ into a bipartite or a tripartite graph with $\chi_{la}=3$, where $n\ge 8$. Let $A=\{v_i\;|\; i \mbox{ is even}\}$ and $B=\{v_i\;|\; i \mbox{ is odd}\}$. \\

\noindent Suppose $n= 4m$. We obtain a partition of $A$ (denoted $\mathscr A$) and of $B$ (denoted $\mathscr B$) such that each block in $\mathscr A$ and in $\mathscr B$ is of size 2. For two distinct vertices $v_i$ and $v_j$ ($i < j$) in each block of $\mathscr A \cup \mathscr B$, we merge them to get a new vertex denoted by $v_{i,j}$ so that the graph obtained is a $4$-regular (not necessarily simple) graph, denoted by $G_{2m}(\mathscr A, \mathscr B)$ (or $G_{2m}$ if no ambiguity). Clearly $G_{2m}$ is bipartite. We preserve the labeling $f$ of $C_{4m}$ to $G_{2m}$. Thus, for $0\le i<j\le 4m-1$,
\[f^+(v_{i,j})=\begin{cases}
6m+4 &  \mbox{ for } i=0;\\
8m+4  &  \mbox{ for } i\ne 0\mbox{ and even};\\
8m+2  &  \mbox{ for } i \mbox{ odd}.
\end{cases}
\]

\noindent Hence, by Lemma~\ref{lem-2part}, we have
\begin{lemma}\label{lem-G4m} For $m\ge 2$ and keeping all notation defined above, $\chi_{la}(G_{2m}(\mathscr A, \mathscr B))=3$.\end{lemma}

\noindent Suppose $n= 4m+2$. Let $\mathscr C=\{v_0, v_{2m+1}\}$ and let $\mathscr A$ and $\mathscr B$ be partitions of $A\setminus\{v_0\}$ and $B\setminus\{v_{2m+1}\}$, respectively, such that each block in $\mathscr A$ and in $\mathscr B$ is of size 2.
For two distinct vertices $v_i$ and $v_j$ ($i < j$) in each block of $\mathscr A \cup \mathscr B\cup\mathscr C$, we merge them to get a new vertex denoted by $v_{i,j}$ so that the graph obtained is a $4$-regular (not necessarily simple) graph, denoted by $G^1_{2m+1}(\mathscr A, \mathscr B,\mathscr C)$ (or $G^1_{2m+1}$ if no ambiguity). Clearly $G^1_{2m+1}$ is tripartite. We preserve the labeling $f$ of $C_{4m+2}$ to $G^1_{2m+1}$. Thus, for $0\le i<j\le 4m+1$,
\[f^+(v_{i,j})=\begin{cases}
6m+6 &  \mbox{ for } i=0, j=2m+1;\\
8m+8  &  \mbox{ for } i\ne 0\mbox{ and even};\\
8m+6  &  \mbox{ for } i \mbox{ odd}.
\end{cases}
\]

\noindent Suppose $n= 4m+1$. Let $\mathscr C=\{v_0\}$ and let $\mathscr A$ and $\mathscr B$ be partitions of $A\setminus\{v_0\}$ and $B$, respectively, such that each block in $\mathscr A$ and in $\mathscr B$ is of size 2.
For two distinct vertices $v_i$ and $v_j$ ($i < j$) in each block of $\mathscr A \cup \mathscr B$, we merge them to get a new vertex denoted by $v_{i,j}$ so that the graph obtained is a (not necessarily simple) graph that has a vertex of degree 2 and each other vertex of degree 4, denoted by $G^2_{2m+1}(\mathscr A, \mathscr B,\mathscr C)$ (or $G^2_{2m+1}$ if no ambiguity). Clearly $G^2_{2m+1}$ is tripartite. We preserve the labeling $f$ of $C_{4m+1}$ to $G^2_{2m+1}$. Thus, $f^+(v_0)=2m+2$ and for $1\le i<j\le 4m$,
\[f^+(v_{i,j})=\begin{cases}
8m+6  & \mbox{ for } i \mbox{ even};\\
8m+4  & \mbox{ for } i \mbox{ odd}.
\end{cases}
\]

\noindent Suppose $n= 4m+3$. Let $\mathscr C=\{v_0, v_{m+1}, v_{3m+2}\}$ and let $\mathscr A$ and $\mathscr B$ be partitions of $A\setminus\{v_0, v_{m+1}, v_{3m+2}\}$ and $B\setminus\{v_{m+1}, v_{3m+2}\}$, respectively, such that each block in $\mathscr A$ and in $\mathscr B$ is of size 2. Note that $m+1\not\equiv 3m+2\pmod{2}$. Merge the 3 vertices in $\mathscr C$ to get a degree 6 vertex $z$.
For two distinct vertices $v_i$ and $v_j$ ($i < j$) in each block of $\mathscr A \cup \mathscr B$, we merge them to get a new vertex denoted by $v_{i,j}$ so that the graph obtained is a (not necessarily simple) graph that has a vertex of degree 6 and each other vertex of degree 4, denoted by $G^3_{2m+1}(\mathscr A, \mathscr B,\mathscr C)$ (or $G^3_{2m+1}$ if no ambiguity). Clearly $G^3_{2m+1}$ is tripartite. We preserve the labeling $f$ of $C_{4m+3}$ to $G^3_{2m+1}$. Thus, $f^+(z)=10m+12$, for $i,j\in [1, 4m+2]\setminus\{m+1, 3m+2\}$ and for $1\le i<j\le 4m+2$,
\[f^+(v_{i,j})=\begin{cases}
8m+10  & \mbox{ for } i\mbox{ even};\\
8m+8  & \mbox{ for } i \mbox{ odd}.
\end{cases}
\]

\noindent Since $\chi(G^s_{2m+1}) = 3$ for $s=1,2,3$, we also have

\begin{lemma}\label{lem-G^s_{2m+1}} For $m\ge 2$ and keeping all notation defined above, $\chi_{la}(G^s_{2m+1}(\mathscr A, \mathscr B, \mathscr C)) = 3$. \end{lemma}

\noindent If we define $\mathscr A$ and $\mathscr B$ suitably, we may obtain a circulant graph.

\noindent{\bf Case (1).}  Consider $n=8k$, $k\ge 2$. Let $\mathscr A=\{(v_{2i}, v_{4k+2i})\;|\; i\in [0, 2k-1]\}$ and\\ $\mathscr B=\{(v_{2j+1},v_{2k+2j+1})\;|\;j\in [0, k-1]\cup [2k, 3k-1]\}$.
In $G_{4k}$, we rewrite the vertex $v_{2i, 4k+2i}$ by $u_{2i}$ for $i\in [0, 2k-1]$; $v_{2j+1, 2k+2j+1}$ by $u_{2j+1}$ for $j\in [0, k-1]$ and by $u_{2j-2k+1}$ for $j\in [2k, 3k-1]$. Thus, the vertex set of $G_{4k}$ is $\{u_i\;|\; i\in [0, 4k-1]\}$.

\noindent Now we are going to look at the neighbors of vertex in $G_{4k}$. First, consider $u_{2j}$, $j\in[0,2k-1]$.
\begin{itemize}

\item $j=0$. The neighbors of $v_{0}$ in $C_{8k}$ are $v_{1}$ and $v_{8k-1}$, and those of $v_{4k}$ are $v_{4k+1}$ and $v_{4k-1}$. So, the neighbors of $u_{0}=v_{0, 4k}$ in $G_{4k}$ are
$v_{1, 2k+1}=u_1$, $v_{6k-1, 8k-1}=u_{4k-1}$, $v_{4k+1, 6k+1}=u_{2k+1}$ and $v_{2k-1,4k-1}=u_{2k-1}$.

\item $j\in[1, k-1]$, the neighbors of $v_{2j}$ in $C_{8k}$ are $v_{2j\pm 1}$ and those of $v_{4k+2i}$ are $v_{4k+2i\pm 1}$. So the neighbors of $u_{2j}=v_{2j, 4k+2j}$ in $G_{4k}$ are $v_{2j+ 1, 2k+2j+1}=u_{2j+1}$, $v_{2j-1, 2k+2j-1}=u_{2j-1}$, $v_{4k+2j+1,6k+2j+1}=u_{2k+2j+1}$ and $v_{4k+2j-1,6k+2j-1}=u_{2k+2j-1}$.

\item $j=k$. The neighbors of $v_{2k}$ in $C_{8k}$ are $v_{2k\pm 1}$ and those of $v_{6k}$ are $v_{6k\pm 1}$. So the neighbors of $u_{2k}=v_{2k, 6k}$ in $G_{4k}$ are  $v_{1,2k+1}=u_{1}$, $v_{2k-1, 4k-1}=u_{2k-1}$, $v_{4k+1,6k+1}=u_{2k+1}$ and $v_{6k-1,8k-1}=u_{4k-1}$.

\item $j\in[k+1, 2k-1]$, the neighbors of $v_{2j}$ in $C_{8k}$ are $v_{2j\pm 1}=v_{2k+2(j-k)\pm 1}$ and those of $v_{4k+2j}$ are $v_{4k+2j\pm 1}=v_{2k+2(k+j)\pm 1}$. So the neighbors of $u_{2j}=v_{2j, 4k+2j}$ in $G_{4k}$ are $v_{2(j-k)+ 1,2k+2(j-k)+ 1}=u_{2j-2k+1}$, $v_{2(j-k)- 1,2k+2(j-k)-1}=u_{2j-2k-1}$, $v_{2(k+j)+1,2k+2(k+j)+1}=u_{2j+1}$ and\break $v_{2(k+j)-1,2k+2(k+j)-1}=u_{2j-1}$.
\end{itemize}

\noindent Next, consider $u_{2i+1}$, $i\in[0,2k-1]$.
\begin{itemize}
\item $i\in[0, k-1]$, the neighbors of $v_{2i+1}$ in $C_{8k}$ are $v_{2i+2}$ and $v_{2i}$, and those of $v_{2k+2i+1}$ are $v_{2k+2i+2}$ and $v_{2k+2i}$.
    So the neighbors of $u_{2i+1}=v_{2i+1, 2k+2i+1}$ in $G_{4k}$ are $v_{2i+2, 4k+2i+2}=u_{2i+2}$, $v_{2i, 4k+2i}=u_{2i}$, $v_{2k+2i+2,6k+2i+2}=u_{2k+2i+2}$ and $v_{2k+2i,6k+2i}=u_{2k+2i}$.

\item $i\in[k, 2k-2]$, the neighbors of $v_{2k+2i+1}$ in $C_{8k}$ are $v_{2k+2i+2}$ and $v_{2k+2i}$, and those of $v_{4k+2i+1}$ are $v_{4k+2i+2}$ and $v_{4k+2i}$.
    So the neighbors of the vertex $u_{2i+1}=v_{2k+2i+1,4k+2i+1}$ in $G_{4k}$ are $v_{2k+2i+ 2,6k+2i+2}=u_{2k+2i+2}$, $v_{2k+2i,6k+2i}=u_{2k+2i}$, $v_{2i+2, 4k+2i+2}=u_{2i+2}$ and $v_{2i,4k+2i}=u_{2i}$.

\item $i=2k-1$, the neighbors of $v_{6k-1}$ in $C_{8k}$ are $v_{6k}$ and $v_{6k-2}$, and those of $v_{8k-1}$ are $v_{0}$ and $v_{8k-2}$.
    So the neighbors of the vertex $u_{4k-1}=v_{6k-1,8k-1}$ in $G_{4k}$ are $v_{2k,6k}=u_{2k}$, $v_{2k-2,6k-2}=u_{2k-2}$, $v_{0, 4k}=u_{0}$ and $v_{4k-2,8k-2}=u_{4k-2}$.
\end{itemize}
\noindent Since $-2k\equiv 2k\pmod{4k}$, we get that
\begin{proposition}\label{pro-4regcir} $G_{4k}(\mathscr A, \mathscr B)\cong C_{4k}(1,2k-1)$, where $\mathscr A$ are $\mathscr B$ are defined above.\end{proposition}

\noindent We now show that when $n$ is sufficiently large, we may repeat the above method to obtain $2^s$-regular graphs that are also circulant. Suppose $n = 2^{2s-1}(t+2)$, $t \ge 0$  and $s\ge 2$. Let $[a,b]_e$ (respectively $[a,b]_d$) denote the set of  consecutive evens (respectively odds) from $a$ to $b$. 

\begin{enumerate}[(a)]

  \item Divide the $n/2=2^{2s-2}(t+2)$ evens in $[0, n-2]_e$ into $2^{2s-2}$ groups of consecutive evens of equal size, which is $t+2$, to get $a_1=[0,2t+2]_e$, $a_2=[2t+4,4t+6]_e$, \ldots, $a_{2^{2s-2}-1}=[n-4t-8,n-2t-6]_e$, $a_{2^{2s-2}}=[n-2t-4,n-2]_e$. Write each $a_k$ as a column. Suppose $S$ is a matrix and $c\in \Z$. For convenience, let $S\oplus c$ be the matrix obtained from $S$ by adding each entry by $c$. So, $a_i=a_1\oplus(i-1)(2t+4)$, $1\le i\le 2^{2s-2}$.

      For iteration 1, take the first 4 groups and form a block matrix $A_1=\begin{pmatrix}a_1 & a_3\\a_2 & a_4\end{pmatrix}$. Actually, $A_1$ is a $(2t+4)\times 2$ matrix.

      For iteration $i\ge 2$, repeat with next $2^{2i-2}$ groups below in the same order of the first $2^{2i-2}$ groups and next $2^{2i-1}$ groups to the right in the same order of the first $2^{2i-1}$ groups, using a total of $2^{2i}$ groups. The resulting matrix is denoted by $A_{i}$. In other word, \[A_{i}=\begin{pmatrix}A_{i-1} & A_{i-1}\oplus 2\times 2^{2i-2}(2t+4)\\A_{i-1}\oplus 2^{2i-2}(2t+4) & A_{i-1}\oplus 3\times 2^{2i-2}(2t+4)\end{pmatrix}.\]
      After $s-1\ge 1$ iteration(s), all the $2^{2s-2}$ groups are occupied. Hence we get a $2^{s-1}(t+2)\times 2^{s-1}$ array $A$.

  \item  Divide the $2^{2s-2}(t+2)$ odds in $[1,n-1]_d$ into $2^{2s-2}$ groups of consecutive odds of equal size, which is $t+2$, to get $b_1=[1,2t+3]_d$, $b_2=[2t+5,4t+7]_d$, \ldots, $b_{2^{2s-2}-1}=[n-4t-7,n-2t-5]_d$, $b_{2^{2s-2}}=[n-2t-3,n-1]_d$. Write each $b_k$ as a row. By using the same format as $A_i$ (change $a$ to $b$, $A$ to $B$), after $s-1\ge 1$ iteration(s), all the $2^{2s-2}$ groups are occupied. Hence we get a $2^{s-1}\times 2^{s-1}(t+2)$ array $B$. Actually, $B_1$ is a $2\times (2t+4)$ matrix.

\item For a fixed $t$, we define a square matrix $\mathscr M_{s-1}$ of order $2^{s-1}(t+2)$ as follows. Let the $(1,2^{s-1}(t+2))$-entry of $\mathscr M_{s-1}$ be 1. If
 row $x$ in $A$ and column $y$ in $B$ has a pair of consecutive integers, then we define the $(x,y)$-entry of $\mathscr M_{s-1}$ by 1, otherwise 0.

\end{enumerate}
We now have the following observations.
\begin{enumerate}[(O1).]
\item When $s=2$. We get that $\mathscr M_{1}$ is\\
$\fontsize{8}{14}\selectfont
\begin{array}{*{2}{c}|*{6}{c}:*{6}{c}|}
& & 1 & 3 &  \cdots & \cdots & 2t+1 & 2t+3 & 4t+9 & 4t+11 & \cdots & \cdots & 6t+9 & 6t+11\\
& & 2t+5 & 2t+7 & \cdots & \cdots & 4t+5 & 4t+7 & 6t+13 & 6t+15 &  \cdots & \cdots & 8t+13 & 8t+15\\\hline
0 & 4t+8 & 1 & 0 & \cdots & \cdots & 0 & 1 & 1 & 0 & \cdots & \cdots & 0 & 1\\
2 & 4t+10 & 1 & 1 & 0 & \cdots & 0 & 0 & 1 & 1 & 0 & \cdots & 0 & 0\\
\vdots & \vdots & 0 & \ddots & \ddots &\ddots & \ddots & \vdots & 0 & \ddots & \ddots &\ddots & \ddots & \vdots\\
\vdots & \vdots & \vdots & \ddots & \ddots &\ddots & \ddots & \vdots & \vdots & \ddots & \ddots &\ddots & \ddots & \vdots\\
2t & 6t+8 & 0 & 0 & \ddots & \ddots & 1 & 0 & 0 & 0 & \ddots & \ddots & 1 & 0\\
2t+2 & 6t+10 & 0 & 0 & \cdots & 0 & 1 & 1 & 0 & 0 & \cdots & 0 & 1 & 1\\\hdashline
2t+4 & 6t+12 & 1 & 0 & \cdots & \cdots & 0 & 1 & 1 & 0 & \cdots & \cdots & 0 & 1\\
2t+6 & 6t+14 & 1 & 1 & 0 & \cdots & 0 & 0 & 1 & 1 & 0 & \cdots & 0 & 0\\
\vdots & \vdots & 0 & \ddots & \ddots &\ddots & \ddots & \vdots & 0 & \ddots & \ddots &\ddots & \ddots & \vdots\\
\vdots & \vdots & \vdots & \ddots & \ddots &\ddots & \ddots & \vdots & \vdots & \ddots & \ddots &\ddots & \ddots & \vdots\\
4t+4 & 8t+12 & 0 & 0 & \ddots & \ddots & 1 & 0 & 0 & 0 & \ddots & \ddots & 1 & 0\\
4t+6 & 8t+14 & 0 & 0 & \cdots &  0 & 1 & 1 & 0 & 0 & \cdots & 0 & 1 & 1\\\hline
\end{array}.$

Clearly the above matrix is circulant. For convenience, we let $\mathscr R$ be the submatrix of $\mathscr M_1$ generated by the first $t+2$ rows and the first $t+2$ columns. Hence, $\mathscr M_1$ is a $2\times 2$ block matrix of the following form:\\
\centerline{$\mathscr M_1=\begin{pmatrix}\mathscr R & \mathscr R \\ \mathscr R  & \mathscr R\end{pmatrix}$.}

%
%
\item Suppose $i\ge 2$. For convenience we denote $2^{2i-2}(2t+4)$ by $m$.\\ Here $A_{i}=\begin{pmatrix}A_{i-1} & A_{i-1}\oplus 2m\\A_{i-1}\oplus m & A_{i-1}\oplus 3m\end{pmatrix}$ and $B_{i}=\begin{pmatrix}B_{i-1} & B_{i-1}\oplus 2m\\B_{i-1}\oplus m & B_{i-1}\oplus 3m\end{pmatrix}$.
We partition $\mathscr M_i$ into 4 submatrices as follows:
\begin{equation}\label{eq-partition}\begin{array}{c:c|c|c}
\multicolumn{2}{c|}{\multirow{2}{1.3cm}{}} & B_{i-1} & B_{i-1}\oplus 2m\\\cdashline{3-4}
\multicolumn{2}{c|}{\multirow{2}{1.3cm}{}} & B_{i-1}\oplus m  & B_{i-1}\oplus 3m\\\hline
A_{i-1} & A_{i-1}\oplus 2m & \mathscr N_1 & \mathscr N_3\\\hdashline
A_{i-1}\oplus m  & A_{i-1}\oplus 3m & \mathscr N_2 & \mathscr N_4
\end{array}. \end{equation}
It is easy to see that all entries of $A_{i-1}\oplus 2m$ are greater than those of $B_{i-1}$. Moreover, the smallest entry of $A_{i-1}\oplus 2m$ is greater than the largest entry of $B_{i-1}\oplus m$ by 1. This pair of numbers corresponds to the top-rightmost entry of $\mathscr N_1$. So $\mathscr N_1=\mathscr M_{i-1}$.

Next $A_{i-1}\oplus 2m$ and $B_{i-1}\oplus 2m$ are the translation of $A_{i-1}$ and $B_{i-1}$ by $2m$. Under our construction we will define the top-rightmost entry of $\mathscr N_3$ by 1. Thus $\mathscr N_3=\mathscr M_{i-1}$.

Similarly, $A_{i-1}\oplus m$ and $B_{i-1}\oplus m$ are the translation of $A_{i-1}$ and $B_{i-1}$ by $m$. Moreover, the smallest entry of $A_{i-1}\oplus m$ is greater than the largest entry of $B_{i-1}$ by 1. Thus the pattern of $\mathscr N_2=\mathscr M_{i-1}$.

Similarly again, $A_{i-1}\oplus 3m$ and $B_{i-1}\oplus 3m$ are the translation of $A_{i-1}$ and $B_{i-1}$ by $3m$. Moreover, the smallest entry of $A_{i-1}\oplus 3m$ is greater than the largest entry of $B_{i-1}\oplus 2m$ by 1. Thus the pattern of $\mathscr N_4=\mathscr M_{i-1}$.

Thus $\mathscr M_i$ is a $2^i\times 2^i$ block matrix whose entries are $\mathscr R$. Clearly it is circulant.

\item We rename the rows and the columns of $\mathscr M_{s-1}$ by $u_0$, $u_2$, $\dots$, $u_{2^{s}(t+2)-2}$ and $u_1$, $u_3$, $\dots$, $u_{2^{s}(t+2)-1}$ in order. By our construction, $\begin{pmatrix} \bigstar & \mathscr M_{s-1}\\\mathscr M_{s-1}^T & \bigstar\end{pmatrix}$ is the adjacency matrix of\break $C_{2^s(t+2)}(1, 2t+3, 2t+5, 4t+7, \ldots, 1+(2^{s-2}-1)(2t+4), 2t+3+(2^{s-2}-1)(2t+4))$ with vertex set $\{u_i\;|\; 0\le i\le 2^{s}(t+2)-1\}$, where $s\ge 2$.

\item Let us define a labeling matrix $\mathscr M$ for the graph corresponding to  $\mathscr M_{s-1}$.
 If row $x$ in $A$ and column $y$ in $B$ has a pair of consecutive integers say, $p$ and $q$, then we define the $(x,y)$-entry of $\mathscr M$ by
   $$m_{x,y}=\begin{cases}p/2 + 1 & \mbox{ if } p=q-1 \\ n - p/2 +1&\mbox{ if }  p=q+1\end{cases}=\begin{cases}(q-1)/2 + 1 & \mbox{ if } q=p+1\\ n - (q-1)/2&\mbox{ if } q=p-1.\end{cases}$$
   Moreover define $m_{1,2^{s-1}(t+2)}=2^{2s-2}(t+2)+1=n/2+1$ (this is in accordance to the $C$-labeling of $C_n$) and $m_{x,y}=*$ for other ordered pairs $(x,y)$.

For each even number $p\in [0, n-2]$, there is an odd number $q=p+1$. Thus $[0, 2^{2s-2}(t+2)]$ are assigned to some entries of $\mathscr M$. Also, for each even number $p\in [2, n-1]$, there is an odd number $q=p-1$. Thus $[2^{2s-2}(t+2)+2, 2^{2s-1}(t+2)]$ are assigned to some entries of $\mathscr M$. Since $m_{1,2^{s-1}(t+2)}=2^{2s-2}(t+2)+1$, all labels in $[1, 2^{2s-1}(t+2)]$ are assigned. Clearly that there are $2^{2s-1}(t+2)$ 1's is $\mathscr M_{s-1}$. So the labeling is bijective.

\item We want to show by induction that {\it if
 row $x$ in $A$ and column $y$ in $B$ has a pair of consecutive integers say, $p$ and $q$, then there is a $y'\ne y$ such that row $x$ in $A$ and column $y'$ in $B$ has a pair of consecutive integers $p$ and $q'$, except the pairs $($row $1$ in $A$, column $1$ in $B)$ and  $($row $1$ in $A$, column $2^{s-1}(t+2)$ in $B)$}. Let $\mathscr R_{a}$ denote the submatrix by taking the first $a$ rows and the first $a$ columns of $\mathscr M$.

It is easy to see (the adjacent 1's in $\mathscr M_1$) that all rows in $\mathscr R_{2(t+2)}$ has the above property except row $t+3$ in $A$. In row $t+3$ of $A$ and column 1 of $B$, there is a pair $p=2t+4$ and $q=2t+5$. We can find that, in row $t+3$ of $A$ and column $t+2$ of $B$, there is a pair $p=2t+4$ and $q'=2t+3$. Next, in row $t+3$ of $A$ and column $t+3$ of $B$, there is a pair $p=6t+12$ and $q=6t+13$. We can find that, in row $t+3$ of $A$ and column $2t+4$ of $B$, there is a pair $p=6t+12$ and $q'=6t+11$.

Suppose $\mathscr R_{2^{i-1}(t+2)}$ has the above property, for $i\ge 2$. By the same argument in (O2), we get the  $\mathscr R_{2^{i}(t+2)}$ also has the above property. By induction, we obtain that $\mathscr M$ has the above property.

Similarly, we can also obtain that {\it if row $x$ in $A$ and column $y$ in $B$ has a pair of consecutive integers say, $p$ and $q$, then there is an $x'\ne x$ such that row $x'$ in $A$ and column $y$ in $B$ has a pair of consecutive integers $p$ and $q'$}.

By the property above, we have that, each column sum of $\mathscr M$ is $2^{s-1}(n+1)$, each row sum of $\mathscr M$ is $2^{s-1}(n+2)$ except the first row sum. Since $m_{1,1}=1$ and  $m_{1,2^{s-1}(t+2)}=n/2+1$, the first row sum of $\mathscr M$ is $2^{s-1}(n+2)-n/2$. Thus, $\mathscr M$ corresponds a local antimagic labeling with color number $3$.
\end{enumerate}

\noindent So we have
\begin{theorem}\label{thm-circulant2} Suppose $s\ge 2$, $t \ge 0$. We have $$\chi_{la}C_{2^s(t+2)}(1, 2t+3, \ldots, 1+(2^{s-2}-1)(2t+4), 2t+3+(2^{s-2}-1)(2t+4))=3.$$
\end{theorem}

\noindent Note that for $n=16$, the new graph obtained by the above method is $K_{4,4}\cong C_8(1,3)$. Although Theorem~\ref{thm-circulant2} is a special case of Theorem~\ref{thm-circulant}, the approach of obtaining this theorem is required in Section~\ref{sec-cycle-union}.

\begin{example} For $n=128$, $s=3$, $t=2$, the labeling matrix $\mathscr M$ is given below.
\[\fontsize{8}{10}\selectfont
\begin{array}{*{4}{>{\raggedleft}p{4mm}}|*{4}{>{\raggedleft}p{4mm}}|*{4}{>{\raggedleft}p{4mm}}|*{4}{>{\raggedleft}p{4mm}}|*{3}{>{\raggedleft}p{4mm}}r|r}
\multicolumn{4}{c|}{\multirow{4}{1.3cm}{$C_{32}(1,7,9,15)$}} & 1 & 3 & 5 & 7 & 17 & 19 & 21 & 23 & 65 & 67 & 69 & 71 & 81 & 83 & 85 & 87 \\
&&& & 9 & 11 & 13 & 15 & 25 & 27 & 29 & 31 & 73 & 75 & 77 & 79 & 89 & 91 & 93 & 95\\
&&& & 33 & 35 & 37 & 39 & 49 & 51 & 53 & 55 & 97 & 99 & 101 & 103 & 113 & 115 & 117 & 119\\
&&& & 41 & 43 & 45 & 47 & 57 & 59 & 61 & 63 & 105 & 107 & 109 & 111 & 121 & 123 & 125 & 127 & \mbox{Sum}\\\hline
0& 16& 64& 80 & 1 & * & * & 121 & 9 & * & * & 97 & 33 & * & * & 89 & 41 & * & * & 65 & 456\\
2 & 18& 66& 82 & 128 & 2 & * & * & 120 & 10 & * & * & 96 & 34 & * & * & 88 & 42 & * & * & 520\\
4& 20& 68& 84 & * & 127 & 3 & * & * & 119 & 11 & * & * & 95 & 35 & * & * & 87 & 43 & * & 520\\
6& 22& 70& 86 & * & * & 126 & 4 & * & * & 118 & 12 & * & * & 94 & 36 & * & * & 86 & 44 & 520\\\cdashline{1-4}\cline{5-21}
8& 24& 72& 88 & 5 & * & * & 125 & 13 & * & * & 117 & 37 & * & * & 93 & 45 & * & * & 85 & 520\\
10& 26& 74& 90 & 124 & 6 & * & * & 116 & 14 & * & * & 92 & 38 & * & * & 84 & 46 & * & * & 520\\
12& 28& 76& 92 &  * & 123 & 7 & * & * & 115 & 15 & * & * & 91 & 39 & * & * & 83 & 47 & * & 520\\
14& 30& 78& 94 & * &  * & 122 & 8 & * & * & 114 & 16 & * & * & 90 & 40 & * & * & 82 & 48 & 520\\\hline
32& 48& 96& 112 & 17 & * &  * & 105 & 25 & * & * & 113 & 49 & * & * & 73 & 57 & * & * & 81 & 520\\
34& 50& 98& 114 & 112 & 18 & * &  * & 104 & 26 & * & * & 80 & 50 & * & * & 72 & 58 & * & * & 520\\
36& 52& 100& 116 & * & 111 & 19 & * &  * & 103 & 27 & * & * & 79 & 51 & * & * & 71 & 59 & * & 520\\
38& 54& 102 & 118 & * & * & 110 & 20 & * &  * & 102 & 28 & * & * & 78 & 52& * & * & 70 & 60 & 520\\\cdashline{1-4}\cline{5-21}
40& 56& 104 & 120 & 21 & * & * & 109 & 29 & * &  * & 101 & 53 & * & * & 77 & 61 & * & * & 69 & 520\\
42& 58& 106& 122 & 108 & 22 & * & * & 100 & 30 & * &  * & 76 & 54 & * & * & 68 & 62 & * & * & 520\\
44& 60& 108& 124 & * & 107 & 23 & * & * & 99 & 31 & * &  * & 75 & 55 & * & * & 67 & 63 & * & 520\\
46& 62& 110& 126 & * & * & 106 & 24 & * & * & 98 & 32 & * &  * & 74 & 56 & * & * & 66 & 64 &  520\\\hline
\multicolumn{4}{r|}{\mbox{Sum}} & 516 & 516 & 516 & 516 & 516 & 516 & 516 & 516 & 516 & 516 & 516 & 516 & 516 & 516 & 516 & 516\\
\end{array}
\]
\qed
\end{example}

\noindent In what follows, let $G-e$ denotes the graph $G$ with an edge $e$ delete. \\

\noindent{\bf Case (2).} Consider $n=8k+4$, $k\ge 2$. Using $\mathscr A=\{(v_{2i}, v_{4k+2i+2})\;|\; i\in [0, 2k]\}$ and $\mathscr B=\{(v_{2j+1},v_{2k+2j+3})\;|\;j\in [0, k]\}\cup \{(v_{4k+2j+5},v_{6k+2j+5})\;|\;j\in [0, k-1]\}$, we obtain a 4-regular bipartite graph, denoted $G_{4k+2}(\mathscr A, \mathscr B)$. Note that when $k=1$ we will get a non-simple graph of order 6. We preserve the labeling $f$ of $C_{8k+4}$ to $G_{4k+2}(\mathscr A, \mathscr B)$. Thus, for $0\le i<j\le 4m+1$,
\[f^+(v_{i,j})=\begin{cases}
12k+10 &  \mbox{ for } i=0;\\
16k+12  &  \mbox{ for } i \mbox{ even};\\
16k+10  &  \mbox{ for } i \mbox{ odd}.
\end{cases}
\]

\noindent By Lemma~\ref{lem-2part}, we conclude that $\chi_{la}(G_{4k+2}(\mathscr A, \mathscr B)) = 3$. It is straightforward to check the conditions of Lemma~\ref{lem-reg}. Together with Lemma~\ref{lem-2part}, we conclude that $\chi_{la}(G_{4k+2}(\mathscr A, \mathscr B)-e)=3$ if $e$ is the edge that receives label 1 or $n$. [Not sure if both resulting graphs are isomorphic.]

\medskip

\noindent{\bf Case (3).} Consider $n=8k+2, k\ge 2$. Using $\mathscr A=\{(v_{2i}, v_{4k+2i})\;|\;i\in[1,2k]\}$,\\ $\mathscr B=\{(v_{2j+1}, $ $v_{4k+2j+3})\;| \;j\in[0,2k-1]\}$ and $\mathscr C=\{(0,4k+1)\}$, we obtain a 4-regular graph with vertices $v_{2,4k+2}, v_{1,4k+3}, v_{0,4k+1}$ form an induced $K_3$ subgraph, denoted $G^1_{4k+1}(\mathscr A, \mathscr B, \mathscr C)$. Preserving the labeling of $C_{8k+2}$ to $G^1_{4k+1}(\mathscr A, \mathscr B, \mathscr C)$, this new graph has $2k$ vertices given by $\mathscr A$ with induced label $16k+8$, $2k$ vertices given by $\mathscr B$ with induced label $16k+6$ and the vertex given by $\mathscr C$ has induced label $12k+6$. Thus, $\chi_{la}(G^1_{4k+1}(\mathscr A, \mathscr B, \mathscr C))=3$. It is straightforward to check the conditions of Lemma~\ref{lem-reg}. Note that $G^1_{4k+1}(\mathscr A, \mathscr B, \mathscr C)-e$ is tripartite if $e$ is the edge that receives label 1 or $n$. Thus, we conclude that $\chi_{la}(G^1_{4k+1}(\mathscr A, \mathscr B, \mathscr C)-e)=3$.

\medskip

\noindent{\bf Case (4).} Consider $n=8k+6, k\ge 2$. Using $\mathscr A=\{(v_{2i}, v_{4k+2+2i})\;|\;i\in[1,2k+1]\}$, $\mathscr B=\{(v_{2j+1},v_{6k+5+2j})\;|\;j\in[0,k]\}\cup\{(v_{2k+3+2j},v_{4k+5+2j})\,|\,j\in[0,k-1]\}$ and $\mathscr C=\{v_0,v_{4k+3}\}$, we obtain a 4-regular graph with vertices $v_{2,4k+2}, v_{1,6k+5}, v_{0,4k+3}$ form an induced $K_3$ subgraph, denoted $G^1_{4k+3}(\mathscr A, \mathscr B, \mathscr C)$. Preserving the labeling of $C_{8k+6}$ to $G^1_{4k+1}(\mathscr A, \mathscr B, \mathscr C)$, this new graph has $2k$ vertices given by $\mathscr A$ with induced label $16k+16$, $2k$ vertices given by $\mathscr B$ with induced label $16k+14$ and the vertex given by $\mathscr C$ has induced label $12k+12$.  Thus, $\chi_{la}(G^1_{4k+3}(\mathscr A, \mathscr B, \mathscr C))=3$. Similar to Case (3), we conclude that $\chi_{la}(G^1_{4k+3}(\mathscr A, \mathscr B, \mathscr C)-e)=3$ if $e$ is the edge that receives label 1 or $n$.

\medskip

\noindent{\bf Case (5).} Consider $n=8k+1, k\ge 2$. Using $\mathscr A=\{(v_{2i}, v_{4k+2i})\;|\; i\in [1, 2k]\}$,\\  $\mathscr B=\{(v_{2j+1},v_{2k+2j+1})\;|\;j\in [0, k-1]\}\cup \{(v_{4k+2j+1},v_{6k+2j+1})\;|\;j\in [0, k-1]\}$ and $\mathscr C=\{v_0\}$, we obtain a tripartite graph denoted $G^2_{4k+1}(\mathscr A, \mathscr B, \mathscr C)$. Preserving the labeling of $C_{8k+1}$ to $G^2_{4k+1}(\mathscr A, \mathscr B, \mathscr C)$, this new graph has $2k$ vertices given by $\mathscr A$ with induced label $16k+6$, $2k$ vertices given by $\mathscr B$ with induced label $16k+4$ and the vertex given by $\mathscr C$ has induced label $4k+2$. Thus, $\chi_{la}(G^2_{4k+1}(\mathscr A, \mathscr B, \mathscr C))=3$. Observe that if $e$ is the edge that receives label 1, then $G^2_{4k+1}(\mathscr A, \mathscr B, \mathscr C)-e$ is a bipartite graph with one pendant having $8k$ edges and partite sets of sizes $2k$ and $2k+1$. By Theorem~\ref{thm-2parteven},  $\chi_{la}(G^2_{4k+1}(\mathscr A, \mathscr B, \mathscr C)-e)\ge 3$. By Lemma~\ref{lem-G-e}, $\chi_{la}(G^2_{4k+1}(\mathscr A, \mathscr B, \mathscr C)-e)\le 3$. Thus, equality holds. Suppose $e$ is the edge that receives label $n$, then  $G^2_{4k+1}(\mathscr A, \mathscr B, \mathscr C)-e$ is a tripartite graph. By Lemmas~\ref{lem-nonreg} and~\ref{lem-G-e}, we conclude that $\chi_{la}(G^2_{4k+1}(\mathscr A, \mathscr B, \mathscr C)-e) =3$.

\medskip

\noindent{\bf Case (6).} Consider $n=8k+5, k\ge 2$. Using $\mathscr A=\{(v_{2i}, v_{4k+2+2i})\;|\;i\in[1,2k+1]\}$, $\mathscr B=\{(v_{2j+1},v_{6k+3+2j})\;|\;j\in[0,k]\}\cup\{(v_{2k+3+2j},v_{4k+3+2j})\,|\,j\in[0,k-1]\}$ and $\mathscr C=\{v_0)\}$, we obtain a tripartite graph denoted $G^2_{4k+3}(\mathscr A, \mathscr B, \mathscr C)$. Preserving the labeling of $C_{8k+5}$ to $G^2_{4k+3}(\mathscr A, \mathscr B, \mathscr C)$, this new graph has $2k+1$ vertices given by $\mathscr A$ with induced label $16k+14$, $2k+1$ vertices given by $\mathscr B$ with induced label $16k+12$ and the vertex given by $\mathscr C$ has induced label $4k+4$. Thus, $\chi_{la}(G^2_{4k+3}(\mathscr A, \mathscr B, \mathscr C))=3$. Similar to Case (5), we conclude that $\chi_{la}(G^2_{4k+3}(\mathscr A, \mathscr B, \mathscr C)-e)=3$ if $e$ is the edge that receives label 1 or $n$.

\medskip

\noindent{\bf Case (7).} Consider $n=8k+3, k\ge 2$. Using $\mathscr A=\{(v_{2i}, v_{4k+2i})\;|\;i\in[1,k]\}\cup\{(v_{2k+2i}, v_{6k+2+2i})\,|$ $i\in[1,k]\}$, $\mathscr B=\{(v_{2j+1},v_{8k+1-2j})\;|\;j\in[0,k-1]\}\cup\{(v_{2k+3+2j},v_{4k+3+2j})\,|\,j\in[0,k-1]\}$ and $\mathscr C=\{v_0, v_{2k+1}, v_{6k+2})\}$,  merge the 3 vertices in $\mathscr C$, we obtain a tripartite graph denoted $G^3_{4k+1}(\mathscr A, \mathscr B, \mathscr C)$. Preserving the labeling of $C_{8k+3}$ to $G^3_{4k+1}(\mathscr A, \mathscr B, \mathscr C)$, this new graph has $2k$ vertices given by $\mathscr A$ with induced label $16k+10$, $2k$ vertices given by $\mathscr B$ with induced label $16k+8$ and the degree 6 vertex $v_{0,2k+1,6k+2}$ has induced label $20k+12$. Thus, $\chi_{la}(G^3_{4k+1}(\mathscr A, \mathscr B, \mathscr C))=3$. Observe that $G^3_{4k+1}(\mathscr A, \mathscr B, \mathscr C)-e$ is a tripartite graph if $e$ is the edge that receives label 1 or $n$. By Lemmas~\ref{lem-nonreg} and~\ref{lem-G-e}, we conclude that $\chi_{la}(G^3_{4k+1}(\mathscr A, \mathscr B, \mathscr C)-e)=3$.

\medskip

\noindent{\bf Case (8).} $n=8k+7, k\ge 2$. Using $\mathscr A=\{(v_{2i}, v_{4k+4+2i})\;|\;i\in[1,k]\}\cup\{(v_{2k+2+2i}, v_{6k+4+2i})\,|\,i\in[1,k+1]\}$, $\mathscr B=\{(v_{4j+1}, v_{4k+3+2j})\;|\;j\in[0,k]\}\cup\{(v_{4j+3},v_{6k+7+2j})\,|\,j\in[0,k-1]\}$ and $\mathscr C=\{v_0, v_{2k+2}, v_{6k+5})\}$,  merge the 3 vertices in $\mathscr C$, we obtain a tripartite graph denoted $G^3_{4k+3}(\mathscr A, \mathscr B, \mathscr C)$. Preserving the labeling of $C_{8k+7}$ to $G^3_{4k+3}(\mathscr A, \mathscr B, \mathscr C)$, this new graph has $2k$ vertices given by $\mathscr A$ with induced label $16k+18$, $2k$ vertices given by $\mathscr B$ with induced label $16k+16$ and the degree 6 vertex $v_{0,2k+2,6k+5}$ has induced label $20k+22$. Thus, $\chi_{la}(G^3_{4k+3}(\mathscr A, \mathscr B, \mathscr C))=3$. Similar to Case (7), we conclude that $\chi_{la}(G^3_{4k+3}(\mathscr A, \mathscr B, \mathscr C)-e)=3$ if $e$ is the edge that receives label 1 or $n$.

\begin{theorem} For $k\ge 2$, if  $e$ is the edge of $G$ that receives label $1$ or $n$, and

\centerline{$\begin{aligned} G\in\{G_{4k+2}(\mathscr A, \mathscr B), & \ G^1_{4k+1}(\mathscr A, \mathscr B, \mathscr C),\ G^1_{4k+3}(\mathscr A, \mathscr B, \mathscr C),\ G^2_{4k+1}(\mathscr A, \mathscr B, \mathscr C),\\ & G^2_{4k+3}(\mathscr A, \mathscr B, \mathscr C),\ G^3_{4k+1}(\mathscr A, \mathscr B, \mathscr C),\ G^3_{4k+3}(\mathscr A, \mathscr B, \mathscr C)\},\end{aligned}$}

\noindent then $\chi_{la}(G)=\chi_{la}(G-e)=3$.
\end{theorem}

\section{Transformation of One-point Union of Cycles}\label{sec-cycle-union}

For $r\ge 2$ and $a_1\ge a_2\ge \cdots\ge a_r\ge 3$, denote by $U(a_1,a_2,\ldots,a_r)$ the one-point union of $r$ distinct cycles of order $a_1, a_2,\ldots,a_r$ respectively. Note that $U(a_1,a_2,\ldots,a_r)$ has $m=a_1+ \cdots + a_r\ge 6$ edges and $m-r+1$ vertices. We shall denote the vertex of maximum degree by $u$, called the {\it central vertex}, and the $2r$ edges incident to $u$ are called the {\it central edges}. For convenience, let $a^{[m]}$ denote a sequence of length $m$ in which all items are $a$, where $m \ge 2$. In~\cite[Theorem 2.5]{LauShiuNg}, the authors completely characterized the local antimagic chromatic number of one-point union of cycles.

\begin{theorem}\label{thm-vertexgluecycles}  Suppose $G=U(a_1,a_2,\ldots,a_r)$. Then $\chi_{la}(G)=  2$ if and only if\\ $G=U((4r-2)^{[r-1]},2r-2)$, $r\ge 3$ or $G=U((2r)^{[\frac{r-1}{2}]}, (2r-2)^{[\frac{r+1}{2}]})$, $r$ is odd. Otherwise, $\chi_{la}(G)= 3$. \end{theorem}

\noindent In this section, we provide an approach to transform a given $U(a_1,a_2,\ldots,a_r)$ into a one-point union of regular bipartite graphs (or bipartite and tripartite graphs) that admits a local antimagic labeling with $\chi_{la}=2$ (or 3).

\subsection{$\chi_{la}(U(a_1,a_2,\ldots,a_r))=2$}\label{sub-chila2}

\noindent First, we consider $G=U((4r-2)^{[r-1]},2r-2), r\ge 3$. Suppose $r\equiv 1$ or $5$ $\pmod{8}, r\ge 9$. Note that $G$ contains an even number of cycles $C_{4r-2}$ and a $C_{2r-2}$. Moreover, in the proof of~\cite[Theorem 2.4]{LauShiuNg}, the local antimagic 2-labeling assigns the edges of the $i$-th copy of $C_{4r-2}$ by $i$, $4r^2-4r+1-i$, $2r-1+i$, $4r^2-6r+2-i$, $4r-2+i$, $4r^2-8r+3-i$, $\ldots$, $4r^2-6r+2+i$, $2r-1-i$ consecutively (for $1\le i\le r-1$) beginning and ending with the two central edges. The edges of the $C_{2r-2}$ are then assigned $2r-1$, $4r^2-6r+2$, $4r-2$, $4r^2-8r+3$, $6r-3$, $4r^2-10r+4$, $\ldots$, $2r^2-3r+1$, $2r^2-r$ consecutively beginning and ending with the two central edges. Note that the central edge labels sum is $4r^2-2r$ while every degree 2 vertex has label $4r^2-4r+1$ and $4r^2-2r$ alternately.\\

\noindent For $s\ge 2$, we shall obtain infinitely many one-point union of $2^s$-regular graphs  (some of which may not be circulant graphs) that has $\chi_{la}=2$ using the local antimagic labeling obtained in Theorem~\ref{thm-vertexgluecycles}. The transformation is done according to the following steps.

\begin{enumerate}[(1)]
\item For each $(i,j)=(1,2), (3,4),  \ldots, (r-2,r-1)$, we perform the following steps:
 \begin{enumerate}[{(1-}1)]
   \item Begin with the $i$-th copy of $C_{4r-2}$, denoted $C^{j/2}$, together with the consecutive edge labels $i$, $4r^2-4r+1-i$, $2r-1+i$, $4r^2-6r+2-i$, $4r-2+i$, $4r^2-8r+3-i$, $\ldots$, $4r^2-6r+2+i$, $2r-1-i$.
   \item Choose an integer $a_{j/2}< 2r-1$ such that $(a_{j/2}, 4r-2)=1$. We add the edges $v_kv_{k+a_{j/2}}$ joining the vertices of $C^{j/2}$ for $0\le k\le 4r-3$ with the indices taken modulo $(4r-2)$ to form another cycle $v_0v_{a_{j/2}}v_{2a_{j/2}}\cdots v_{(4r-3)a_{j/2}}v_0$ of order $4r-2$.
   \item Label the edges of $v_0v_{a_{j/2}}v_{2a_{j/2}}\cdots v_{(4r-3)a_{j/2}}v_0$ by $j$, $4r^2-4r+1-j$, $2r-1+j$, $4r^2-6r+2-j$, $4r-2+j$, $4r^2-8r+3-j$, $\ldots$, $4r^2-6r+2+j,2r-1-j$ consecutively. This gives us an edge labeled $C_{4r-2}(1,a_{j/2})$. Note that for $j\ne j'$, $a_{j/2}$ and $a_{j'/2}$ may be equal.
 \end{enumerate}

\item Consider the $C_{2r-2}$ with consecutive edge labels $2r-1$, $4r^2-6r+2$, $4r-2$, $4r^2-8r+3$, $6r-3$, $4r^2-10r+4$, $\ldots$, $2r^2-3r+1$, $2r^2-r$. Since $2r-2\equiv 0\pmod{8}$, we now transform $C_{2r-2}$ according to steps in Section~\ref{cycles} (that give the graph $G_{2m}$ in Lemma~\ref{lem-G4m}, or the graph in Case (1)) to obtain a simple graph of order $r-1$, denoted $G_{r-1}$, which may or may not be a circulant graph.

\item Identify the vertex $v_0$ of each of the $\frac{r-1}{2}$ copies of $C_{4r-2}(1,a_h), h=1,2,\ldots,\frac{r-1}{2},$ and the vertex with subscript $0$, say $x$, of $G_{r-1}$ and name this merged vertex by $v_0$. Denote the new graph obtained by $G((4r-2)^{[\frac{r-1}{2}]},r-1)$.
\item Label the vertices with their incident edge labels sum.
\end{enumerate}

\noindent We have the following observations.

\begin{enumerate}[(O1)]
  \item The graph $C_{4r-2}(1,a_{j/2})$ has $v_0$ with label $2(2r-1)$ while vertices $v_1$ to $v_{4r-3}$ have labels $2(4r^2-4r+1)$ and $2(4r^2-2r)$ alternately.
  \item The graph $G_{r-1}$ has $x$ with label $(2r-1) + (2r^2-r) + (4r^2-2r) = 6r^2-r-1$ while other consecutive vertices have labels $2(4r^2-4r+1)$ and $2(4r^2-2r)$ alternately.
  \item In $G((4r-2)^{[\frac{r-1}{2}]},r-1)$, vertex $v_0$ and all other vertices of even distance away have vertex label $2(4r^2-2r)$ while the remaining vertices have label $2(4r^2-4r+1)$.
\end{enumerate}

\noindent Consequently, we have obtained a one-point union of $\frac{r-1}{2} + 1$ copies of 4-regular bipartite graphs (possibly with exactly one copy is non-circulant) with $\chi_{la} = 2$.

\begin{rem}\label{rem-2^s} As observed in Section~\ref{cycles} Case (1), if $r = 2^{2s-2}(t+2)+1, t\ge 0$, we can transform the cycle $C_{2r-2}$ repeatedly to obtain a labeled $2^s$-regular bipartite circulant for $s\ge 3$. Consequently, we shall also take $2^{s-1}(t+2)$ groups of $2^{s-1}$ copies of $C_{4r-2}$ and form a $2^s$-regular bipartite circulant according to the approach in Section~\ref{sec-cir}. This new graph will have a vertex $z$ of degree $2r-2+2^s$ and all other vertices of degree $2^s$. Note that all the original central edge labels are now adjacent to $z$ to still have total sum of $4r^2-2r$ while the other $2^s$ edge labels adjacent to $z$ (necessarily from the $C_{2r-2}$) contributed a total sum of $(2^{s-1}-1)(4r^2-2r)$. So, vertex $z$ and all other vertices of even distance away have label $2^{s-1}(4r^2-2r)$ while the remaining vertices have label $2^{s-1}(4r^2-4r+1)$. Thus, the above approach allows us to obtain many non-isomorphic bipartite graphs, each of which is a one-point union of $2^s$-regular bipartite circulants, with $\chi_{la}=2$.
\end{rem}

\begin{rem}\label{rem-37mod8} For $r\equiv 3, 7\pmod{8}$, $2r-2\equiv 4\pmod{8}$. Similar to Step~(1), we transform every pair of the $r-1$ copies of $C_{4r-2}$ to a $C_{2r-1}(1,a_h)$ for $h=1,2,\ldots,(r-1)/2$, $1< a_h < 2r-1$ and $(a_h,4r-2)=1$. Similar to Step~(2), we transform the $C_{2r-2}$ according to Section~\ref{cycles} Case~(2). Finally, by approaches similar to Steps~(3) and (4), we obtain various one-point union of $(r-1)/2$ copies of 4-regular bipartite circulants of order $2r-1$ and a 4-regular bipartite graph of order $r-1$ with $\chi_{la}=2$. The vertex $x$ of degree $2r+2$ and all vertices of even distance away have label $2(4r^2-2r)$ while the remaining vertices have label $2(4r^2-4r+1)$.
\end{rem}

\begin{example} For $r\equiv1,5\pmod{8}$, we can have $U(130^{[32]},64)$ with $r=33$. If $s=2$, we get a one-point union of 16 copies of 4-regular circulants $C_{130}(1,a_h)$ ($1< a_h < 65$ with $(a_h,130)=1$ for $h=1,2,\ldots,16$) and a 4-regular circulant $C_{32}(1,15)$. If $s=3$, we get a one-point union of 8 copies of 8-regular circulants $C_{130}(1,a_i,a_j,a_k)$ ($a_i, a_j$ and $a_k$ are distinct with $1 < a_i, a_j, a_k < 65$ and $(a_i,130) = (a_j,130) = (a_k,130)=1$) and an 8-regular circulant $C_{16}(1,3,5,7)= K_{8,8}$. The resulting graph is of local antimagic chromatic number 2. \\

\noindent For $r\equiv3,7\pmod{8}$, we can have $U(58^{[14]},28)$ with $r=15$, $1 < a_h < 29$ and $(a_h,58)=1$ to get a one-point union of 7 copies of 4-regular circulants $C_{58}(1,a_h)$ and a 4-regular bipartite graphs of order 14. The resulting graph is of local antimagic chromatic number 2. \qed
\end{example}

\noindent Next, consider $U((2r)^{[\frac{r-1}{2}]}, (2r-2)^{[\frac{r+1}{2}]})$, $r$ odd. In the proof of~\cite[Theorem 2.4]{LauShiuNg}, the local antimagic 2-labeling assigns the edges of the $i$-th copy ($1\le i\le (r-1)/2$) of $C_{2r}$ by $i$, $2r^2-r-i$, $2r+i$, $2r^2-3r-i$, $4r+i$, $2r^2-5r-i$, $\ldots$, $r^2-r+i$, $r^2-i$, $\ldots$, $2r^2-6r-i$, $5r-i$, $2r^2-4r+i$, $3r-i$, $2r^2-2r+i$, $r-i$, while the edges of the $k$-th copy ($k=(r+1)/2 + j, 0\le j\le (r-1)/2$) of $C_{2r-2}$ are assigned $r+j$, $2r^2-2r-j$, $3r+j$, $2r^2-4r-j$, $5r+j$, $2r^2-6r-j$, $\ldots$, $r^2+r-j$, $r^2+j$, $\ldots$, $2r^2-7r+j$, $6r-j$, $2r^2-5r+j$, $4r-j$, $2r^2-3r+j$, $2r-j$. Note that the central edge labels sum is $2r^2+r$ while every degree 2 vertex has label $2r^2-r$ and $2r^2+r$ alternately. \\

\noindent Suppose $r\equiv 1,5\pmod{8}, r\ge 9$. We shall apply Step (1) (in Subsection~\ref{sub-chila2}) to the $\frac{r-1}{2}$ copies of $C_{2r}$ to obtain $\frac{r-1}{4}$ copies of bipartite circulants, namely $C_{2r}(1,a_{h})$ ($1 < a_h < r$ with $(a_h,2r)=1$ for $h=1,2,\ldots, (r-1)/4$). For the  $\frac{r-1}{2}$ copies of $C_{2r-2}$, we also apply Step (1) to obtain $\frac{r-1}{4}$ copies of bipartite circulants, namely $C_{2r-2}(1,b_h)$ ($1 < b_h < r-1$ with $(b_h,2r-2)=1$ for $h=1,2,\ldots, (r-1)/4$). Next, we apply Step (2) to the remaining copy of $C_{2r-2}$ to obtain a 4-regular graph $G_{r-1}$. Finally, apply Steps (3) and (4). Now, vertex $x$ and all other vertices of even distance away have label $2(2r^2+r)$ while all other vertices have label $2(2r^2-r)$. Thus, we have obtained various one-point union of $\frac{r-1}{2}+1$ copies of 4-regular bipartite graphs (possibly with exactly one copy is non-circulant) with $\chi_{la}=2$.

\begin{rem}\label{rem-2^s-v2} Similar to Remark~\ref{rem-2^s}, if $r=2^{2s-2}(t+2)+1,t\ge 0$, we can transform a cycle $C_{2r-2}$ repeatedly to obtained a labeled $2^s$-regular bipartite circulant for $s\ge 2$. Consequently, we shall also take $2^{s-2}(t+2)$ groups of $2^{s-1}$ copies of $C_{2r}$ and remaining $C_{2r-2}$ and form a $2^s$-regular bipartite circulant according to the approach in Section~\ref{sec-cir}. This new graph also has a vertex $z$ of degree $2r-2+2^s$ and all other vertices of degree $2^s$. Moreover, vertex $z$ and all other vertices of even distance away have label $2^{s-1}(2r^2+r)$ while all other vertices have label $2^{s-1}(2r^2-r)$. Consequently, we obtain various non-isomorphic bipartite graphs, each of which is a one-point union of $2^s$-regular bipartite circulant graphs, with $\chi_{la}=2$.
\end{rem}

\begin{rem}\label{rem-37mod8-v2} Similar to Remark~\ref{rem-37mod8}, when $r\equiv 3,7\pmod{8}$, we can obtain various one-point union of $(r-1)/4$ copies of 4-regular bipartite circulants of order $r$, $(r-1)/4$ copies of 4-regular bipartite circulants of order $r-1$ and a 4-regular bipartite graph of order $r-1$ with $\chi_{la}=2$. The vertex $x$ of degree $2r+2$ and all all other vertices of even distance away have label $2(2r^2+r)$ while all other vertices have label $2(2r^2-r)$.
\end{rem}

\begin{example} Using $U(34^{[8]},32^{[9]})$ with $r=17$, $s=2$, we get a one-point union of 4 copies of $C_{34}(1,a_h)$, $4$ copies of $C_{32}(1,a_{h'})$ $(a_h, a_{h'}\in\{3,5,7,11,13,15\})$ and a copy of 4-regular circulant $C_{16}(1,7)$. The resulting graph is of local antimagic chromatic number 2. \qed
\end{example}

\noindent Consider $G=U((4r-2)^{[r-1]},2r-2)$ for $r\equiv 0\pmod{2}, r\ge 8$. Since $4r-2\equiv 6\pmod{8}$, we transform each copy of $C_{4r-2}$ as in Section~\ref{cycles} Case (4) to get a 4-regular tripartite graph. If $2r-2 \equiv 2\pmod{8}$, we transform $C_{2r-2}$ as in Section~\ref{cycles} Case (3) to get a 4-regular tripartite graph. If $2r-2 \equiv 6\pmod{8}$, we also transform $C_{2r-2}$ as in Section~\ref{cycles} Case (4) to get a 4-regular tripartite graph. Merging all the vertices with $0$ in its subscript, we obtain a one-point union of 4-regular tripartite graphs that has $\chi_{la} = 3$. \\

\noindent Finally, suppose $U((2r)^{[(r-1)/2]}, (2r-2)^{[(r+1)/2]})$ for $r\equiv 3,7\pmod{8}, r\ge 9$. Since $2r\equiv 6\pmod{8}$, we transform each $C_{2r}$ as in Section~\ref{cycles} Case (4) to get a 4-regular tripartite graph. Since $2r-2\equiv 4\pmod{8}$, we transform each $C_{2r-2}$ as in Section~\ref{cycles} Case (2) to get a 4-regular bipartite graph. Merging all the vertices with $0$ in its subscript, we obtain a one-point union of 4-regular bipartite and tripartite graphs that has $\chi_{la} = 3$. \\

\noindent 
Our next remark is for all graphs above including those 3-partite graphs in the 2 paragraphs above.
\begin{rem} For each graph $G$ obtained above, if $f(e)=1$ or $n$, then Theorem~\ref{thm-2part}, Lemmas~\ref{lem-nonreg} and~\ref{lem-G-e} together imply that $\chi_{la}(G-e)=3$.  \end{rem}

\subsection{$\chi_{la}(U(a_1,a_2,\ldots,a_r))=3$}\label{sub-chila3}

\noindent We only consider $a_i\ge 16$ for each $i$, $1\le i\le r$ and $r\ge 2$.
Denote the consecutive edges of subgraph $C_{a_i}$ by $e_{s_i+1}, e_{s_i+2}, \ldots, e_{s_i+a_i}$ such that $s_1=0$, $s_i=a_1+a_2+\cdots+a_{i-1}$ for $i\ge 2$. Moreover, for $i\ge 1$, $e_{s_i+1}$ and $e_{s_i+a_i}$ are the central edges of $C_{a_i}$. In the proof of~\cite[Theorem 2.4]{LauShiuNg}, a required local antimagic 3-labeling is given by $f: E(G)\to [1,m]$ $(m = \sum^r_{i=1} a_i\ge 32$) such that
\begin{enumerate}[(1)]
  \item $f(e_i) = i/2$ for even $i$,
  \item $f(e_i) = m - (i-1)/2$ for odd $i$
\end{enumerate}
with the degree $2r$ vertex $u$ and every two adjacent degree 2 vertices have labels $m+1$ and $m$ respectively. Moreover, $f^+(u)\ge f(e_1)+[f(e_{a_1})+f(e_{a_1+1})]+f(e_m)\ge m +m+m/2\ge 2m+16$.

\noindent Consider the following two possibilities.
\begin{enumerate}[(a)]
  \item Each $a_i\equiv 0,4\pmod{8}$. Without loss of generality, let $G=U((8k)^{[n_1]},(8k+4)^{[n_2]})$ for $k\ge 2$ and $r=n_1+n_2\ge2$ and $m=8kr+4n_2$. We now transform each $C_{8k}$ and $C_{8k+4}$ according to Section~\ref{sec-cir} Cases (1) or (2) accordingly by preserving the labeling $f$ above to get a one-point union of bipartite graphs. Thus, this graph  has a degree $4r$ vertex with label $> 2m+16$ and each other degree 4 vertex has label $2m+2$ and $2m$ alternately. Therefore, the resulting graph $H$ is a bipartite graph with $\chi_{la}(H)\le 3$. Note that $H$ has bipartition $(V_1, V_2)$ with $|V_1| = 2kr+n_2$ and $|V_2| = 2kr+n_2-r+1$. Observe that $\frac{m}{|V_2|} = 4 + \frac{4(r-1)}{2kr+n_2-r+1}$, not an integer if $k\ge 3$. Suppose $k=2$, then $\frac{m}{|V_2|} = 5+\frac{r-n_2+5}{3r+n_2-1}$, also not an integer. By Theorem~\ref{thm-2part}, $\chi_{la}(H)\ne 2$. Thus, $\chi_{la}(H)=3$.

  \item At least one of $a_i\not\equiv 0,4\pmod{8}$.  We now transform each of $C_{a_i}$ according to Section~\ref{sec-cir} Cases (1) to (8) accordingly by preserving the labeling $f$ above to get a one-point union of tripartite (and possibly bipartite) graphs. The unique vertex that is incident to all the central edges has degree at least $2r$ with label greater than $2m+16$ and each other adjacent degree  4 vertices have labels $2m+2$ and $2m$, respectively. Therefore, the resulting graph $H$ is a tripartite graph with $\chi_{la}(H)= 3$.
\end{enumerate}

\begin{rem} Note that each graph in Part (a) with an edge deleted is bipartite with size $m=8kr+4n_2-1$ and a partite set of size $2kr+n_2$. Since $(8k4+4n_2-1)/ (2kr+n_2)$ is not an integer, by Theorem~\ref{thm-2part}, the graph obtained has $\chi_{la}\ge 3$. By Lemmas~\ref{lem-nonreg} and~\ref{lem-G-e}, each graph $G$ in Part (a) has $\chi_{la}(G-e)=3$ if $f(e) =1$ or $n$. Moreover, by Lemmas~\ref{lem-nonreg} and~\ref{lem-G-e}, each graph $G$ in Part (b) also has $\chi_{la}(G-e)=3$ if $f(e)=1$ or $n$. \end{rem}

\begin{theorem} There are infinitely many one-point union of regular graphs (possibly one copy is not circulant) with $\chi_{la}=2$. \end{theorem}

\begin{theorem} There are infinitely many one-point union $4$-regular bipartite (possibly with tripartite) graphs (with at most one edge deleted) having $\chi_{la}=3$.  \end{theorem}

\section{Concluding Remarks and Open Problems}

In this paper, we first give a sufficient condition for a graph with one pendant with $\chi_{la}\ge 3$. A necessary and sufficient condition for a graph to have $\chi_{la}=2$ is thus obtained. We then obtained infinitely many bipartite circulants with $\chi_{la}=3$.

\begin{question}\label{que-nocoprime}
What is the local antimagic chromatic number of $C_{2n}(a_0, a_1,  \ldots, a_t)$ if all $a_j$ are odd but not each of them is coprime with $2n$?
\end{question}

\begin{question}\label{que-odd}
What is the local antimagic chromatic number of $C_{m}(1, a_1,  \ldots, a_t)$ when $m$ is odd?
\end{question}

\noindent By transforming an $n$-cycle, $n\ge 16$, we obtained infinitely many bipartite and tripartite graphs with $\chi_{la} = 2,3$. Only eight different transformations are given. However, as $n$ becomes arbitrarily large, it is likely that there are many different transformation of cycles that would give infinitely many bipartite and tripartite graphs with $\chi_{la}=3$. Let $\mathscr G$ be the set of all the graphs that can be obtained through all possible transformations of cycles.

\begin{problem} Determine the local antimagic chromatic number of all the graphs in $\mathscr G$. \end{problem}

\noindent Applying the transformation of cycles to one-point union of cycles, we then obtained infinitely many one-point union of regular bipartite graphs (possibly with all except one non-circulant) with $\chi_{la}=2$ and infinitely many one-point union of tripartite (possibly with bipartite) graphs with $\chi_{la}=3$. Note that all the published results on 2-connected bipartite graphs have $\chi_{la}=2$ or 3.
%
Moreover, all the published results on tripartite graphs with at most one pendant have $\chi_{la}=3$ or 4 (see~\cite{Arumugam,LNS-CM,LNS,LauNgShiu-chila,LauShiuNg}).  

\begin{problem} Characterize biparite graphs with  $\chi_{la} = 2$ or 3. \end{problem}

\begin{problem} Characterize tripartite graphs with  $\chi_{la} = 3$ or 4. \end{problem}

\begin{question} Does there exist a bipartite or a tripartite graph with relatively small number of pendants to have arbitrarily large local antimagic chromatic number? \end{question}

\


\end{document}